\documentclass[sn-mathphys-num]{sn-jnl}

\usepackage[T1]{fontenc}																
\usepackage[utf8]{inputenc}															
\usepackage[english]{babel}															
\usepackage{algorithm}%
\usepackage{algorithmicx}%
\usepackage{algpseudocode}%
\usepackage{amsfonts}																
\usepackage{amsmath}																
\usepackage{amssymb}																
\usepackage{amsthm}																
\usepackage[title]{appendix}%
\usepackage{bm}																	
\usepackage{bbm}																	
\usepackage{booktabs}
\usepackage{braket}																	
\usepackage[font=small,format=hang,labelfont={sf,bf}]{caption}									
\usepackage{changepage}															
\usepackage{comment}																
\usepackage{empheq}																
\usepackage{enumitem}																
\usepackage{fancyhdr}																
\usepackage[flushmargin,hang,symbol]{footmisc} 
\usepackage{geometry}	
\usepackage{guit}																	
\usepackage{indentfirst}																
\usepackage{lipsum}																	
\usepackage{listings}																
\usepackage{manyfoot}%
\usepackage{mathrsfs}																
\usepackage{mathtools}																
\usepackage{multirow}
\usepackage{syntonly}
\usepackage{textcomp}
\usepackage[title]{appendix}
\usepackage{url}
\usepackage[english]{varioref}															
\usepackage{xcolor}																	
\usepackage{hyperref}																
\hypersetup{
colorlinks,
citecolor=black,
filecolor=black,
linkcolor=black,
urlcolor=black
}

\newcommand{\N}{\mathbb{N}} 

\newcommand{\R}{\mathbb{R}}

\newcommand{\F}{\mathbb{F}}

\newcommand{\BB}{\mathcal{B}}
\newcommand{\CC}{\mathcal{C}}

\newcommand{\FF}{\mathcal{F}}

\newcommand{\HH}{\mathcal{H}}

\newcommand{\LL}{\mathcal{L}}

\newcommand{\NN}{\mathcal{N}}

\newcommand{\UU}{\mathcal{U}}

\newcommand{\UUU}{\mathscr{U}}

\renewcommand{\epsilon}{\varepsilon} 
\renewcommand{\theta}{\vartheta}
\renewcommand{\rho}{\varrho}
\renewcommand{\phi}{\varphi}

\newcommand{\T}{\mathsf{T}} 

\renewcommand{\u}{\bm{u}} 

\renewcommand{\d}{\textup{d}} 
\newcommand{\D}{\textup{D}} 

\newcommand{\1}{\mathbbm{1}} 

\newcommand{\eqnot}{\coloneqq} 

\DeclarePairedDelimiter{\rounds}{{\mathopen{(}}}{{\mathclose{)}}} 
\DeclarePairedDelimiter{\bracks}{\lbrack}{\rbrack} 
\DeclarePairedDelimiter{\braces}{\lbrace}{\rbrace} 
\DeclarePairedDelimiter{\abs}{\lvert}{\rvert} 
\DeclarePairedDelimiter{\norma}{\lVert}{\rVert} 

\newcommand{\LLped}[3]{\LL_{#3}^{#1}\rounds{#2}} 

\newcommand{\xto}[2]{\xrightarrow[#2]{#1}} 

\DeclareMathOperator{\probP}{\mathbf{P}} 
\DeclareMathOperator{\E}{\mathbf{E}} 

\DeclareMathOperator{\meas}{\mathbf{m}} 
\DeclareMathOperator{\dist}{\mathbf{d}} 

\DeclareMathOperator{\intpart}{int} 
\DeclareMathOperator{\sgn}{sgn} 

\DeclareMathOperator*{\argmin}{arg\,min} 

\newenvironment{sistema}
{\left\lbrace\begin{array}{@{}l@{}}}
{\end{array}\right.}	

\theoremstyle{thmstyleone}%
\newtheorem{theorem}{Theorem}
\newtheorem{corollary}[theorem]{Corollary}
\newtheorem{lemma}[theorem]{Lemma}%

\theoremstyle{thmstyletwo}%
\newtheorem{example}{Example}%
\newtheorem{remark}{Remark}%

\theoremstyle{thmstylethree}%
\newtheorem{assumption}{Assumption}
\newtheorem{definition}{Definition}%
\newtheorem{notation}{Notation}
\newtheorem{problem}{Problem}

\begin{document}

\title[Equilibrium strategies for constrained time\:\!-inconsistent control problems]{\textcolor{black}{Equilibrium strategies for constrained time\:\!-inconsistent control problems}} 


\author[1]{\fnm{Elisa} \sur{Mastrogiacomo}}\email{elisa.mastrogiacomo@uninsubria.it}
\equalcont{These authors contributed equally to this work.}

\author*[2]{\fnm{Marco} \sur{Tarsia}}\email{marco.tarsia@uninsubria.it}
\equalcont{These authors contributed equally to this work.}


\affil[1,2]{\orgdiv{Dipartimento di Economia (DiECO)}, \orgname{Università degli Studi dell’Insubria}, \orgaddress{\street{Via Monte Generoso 71}, \city{Varese}, \postcode{21100}, \country{Italy}}}




\abstract{\textcolor{black}{This paper addresses the issue of time inconsistency in recursive stochastic control problems, in which the forward state process evolves under the influence of an additional recursive utility system. Through an adaptation of Ekeland's variational principle, we establish necessary conditions for subgame\:\!-perfect (Nash) equilibrium strategies, formulated in terms of a Hamiltonian framework defined via coupled backward stochastic differential equations. To illustrate the scope of the results, we consider a constrained portfolio management problem with a finite deterministic horizon and non-exponential discounting, demonstrating the applicability of the proposed methodology in a financial context. The class of admissible constraints examined includes, in particular, the imposition of risk limits on the terminal wealth.}}





\maketitle


\section{Introduction}

\textcolor{black}{The primary objective of this paper is to derive necessary conditions ensuring the existence of equilibrium strategies of closed\:\!-\:\!loop (or feedback) type in time\:\!-inconsistent recursive stochastic control problems subject to suitable constraints, with particular attention to their potential applications in financial decision models. In this vein, the present work extends the unconstrained analysis of~\cite{mastrogiacomotarsia23} by incorporating state\:\!-dependent constraints through an additional recursive utility, thereby introducing a new layer of coupling between the control and the utility dynamics.}

\textcolor{black}{A closed\:\!-\:\!loop equilibrium control can be formally described as a state\:\!-feedback mapping that assigns, at each time, an admissible control action depending on the current state of the system. This structure allows the controller to adjust decisions dynamically as the state evolves, independently of the initial condition, and is particularly advantageous in uncertain or fluctuating environments. By contrast, an open\:\!-\:\!loop equilibrium strategy is determined once and for all at the initial time, depending solely on the starting state, and does not adapt to subsequent changes in the system. We refer the reader to~\cite{weiyongyu17} and~\cite{hamaguchi21} for a detailed comparative discussion.}

\textcolor{black}{In this perspective, time inconsistency refers to the failure of Bellman's principle of optimality and, consequently, of the standard dynamic programming framework. Heuristically, it arises whenever the preferences or objective of a decision-maker evolve over time, so that an action deemed optimal at a given moment may no longer be regarded as optimal from a later standpoint. The resulting lack of intertemporal alignment prevents the problem from being solved through classical dynamic programming techniques, since the optimal strategy cannot be preserved across different time instants. This phenomenon has been extensively investigated in the foundational works of Ekeland and coauthors,~\cite{ekelandlazrak06, ekelandmbodjipirvu12, ekelandpirvu08}, where the focus lies on unconstrained time\:\!-inconsistent stochastic control problems.}

\textcolor{black}{In those contributions, the notion of optimality is replaced by the concept of subgame\:\!-perfect equilibrium, while the stochastic Pontryagin-type maximum principle---in the sense of the seminal work~\cite{yongzhou99}---provides an analytical foundation ensuring time consistency in the classical linear-quadratic and Merton-type portfolio settings, leading to characterizations of possibly non-unique equilibria. This research direction has since evolved into a well-developed literature encompassing both stochastic and deterministic formulations of time\:\!-inconsistent control. We refer, among others, to~\cite{merton71, hujinzhou12, bjorkmurgoci14, bjorkmurgocizhou14, bjorkkhapkomurgoci17, hujinzhou17, yong14} for the stochastic case, and to~\cite{strotz55, pollak68, ekelandlazrak08, ekelandlazrak10} for deterministic counterparts.}

\textcolor{black}{A substantial body of research has been devoted to time\:\!-inconsistent control problems and to the derivation of necessary conditions ensuring the existence of equilibrium strategies; see, for instance,~\cite{weiyongyu17}. These studies are of significant importance, since many practical models in optimization and control inherently exhibit time inconsistency. In the present work, we extend these principles by establishing a generalized maximum principle for a constrained recursive utility optimization problem with non-exponential discounting. Our contribution advances the existing literature in two main directions: first, by incorporating recursive utilities into the time\:\!-inconsistent framework; and second, by explicitly accounting for constraints.}

\textcolor{black}{The analysis is framed in terms of a class of evolutionary equations, either forward or backward in time, leading to performance functionals that depend on fixed initial parameters. This formulation is conceptually analogous, for instance, to that adopted in~\cite{weiyongyu17}. A markedly different yet complementary viewpoint is offered in~\cite{hamaguchi21}, where the analysis is developed within the broader setting of backward stochastic Volterra integral equations. The resolution of such equations requires studying processes known as diagonal processes, allowing both temporal arguments to vary simultaneously. Although this approach is more general and extends beyond the techniques employed here, to the best of our knowledge the present study represents the first attempt to address the constrained case within a time\:\!-inconsistent framework of this type.} 

\textcolor{black}{A complementary line of research has recently developed an extended dynamic programming principle for sophisticated agents in non\:\!-Markovian settings. 
In particular,~\cite{hernandezpossamai23} introduces a game-theoretic framework leading to a rigorous dynamic programming formulation and a characterization of equilibria through an infinite family of backward stochastic differential equations (BSDEs), in analogy to the classical Hamilton\:\!--Jacobi\:\!--Bellman approach, and is further extended in~\cite{hernandezpossamai24} to a finite\:\!-horizon moral hazard problem with a time\:\!-inconsistent agent, unveiling a novel class of control problems with Volterra\:\!-type features and infinite\:\!-dimensional constraints. These works provide a dynamic-programming perspective distinct from the maximum-principle approach pursued in the present study.}

\textcolor{black}{The paper focuses on the constrained case, which constitutes the most technically demanding component of the investigation, and builds upon the unconstrained framework developed in~\cite{ekelandpirvu08} and, more recently, in~\cite{mastrogiacomotarsia23}. In particular,~\cite{mastrogiacomotarsia23} investigates equilibrium strategies for recursive stochastic control problems in the absence of additional constraints, extending the classical spike variation technique to a recursive setting. The resulting formulation provides both necessary and sufficient equilibrium conditions expressed through a second-order Hamiltonian function defined via backward stochastic differential equations, incorporating first- and second-order adjoint equations. For a comprehensive overview of the existing literature on this topic, we refer to~\cite{mastrogiacomotarsia23} and the references therein, including~\cite{ekelandlazrak06},~\cite{ekelandpirvu08}, and~\cite{hu17} (see also~\cite{duffieepstein92},~\cite{elkarouipengquenez01}, and~\cite{elkarouipengquenez97}).}

\textcolor{black}{The constrained formulation under consideration arises from the introduction of an additional utility, leading to an auxiliary recursive system that imposes a state\:\!-dependent restriction on the control process. To tackle this problem, we adapt Ekeland's variational principle, discussed in~\cite{ekeland74}, within a penalization framework. Specifically, we first apply the analytical techniques developed in~\cite{mastrogiacomotarsia23} to derive necessary equilibrium conditions for a sequence of unconstrained approximating problems. The corresponding conditions for the original constrained problem are then recovered by passing to the limit under appropriate boundedness assumptions. An essential technical step in this construction involves the use of the distance function associated with a closed subset of the Euclidean space.}

\textcolor{black}{In this setting, we also establish a transversality-type condition. Owing to the distinctive structure of the framework, no definitive information can be obtained regarding the sign of the scalar multiplier. As a consequence, the classical expression of the Hamiltonian function cannot be directly retrieved within this constrained framework. Although this limitation formally departs from the classical form one would expect for the necessary conditions, it does not reduce their analytical content, as the equilibrium system remains fully characterized through the adjoint equations and the generalized Hamiltonian inequality.}

\textcolor{black}{In contrast, such indeterminacy is absent in the unconstrained framework explored in~\cite{mastrogiacomotarsia23}, where no Lagrange multiplier emerges and the transversality condition retains its classical form. The introduction of a genuine scalar multiplier in the present formulation thus marks a structural departure from that earlier setting.}

\textcolor{black}{Rather, it prevents a straightforward reduction to the standard Lagrange\:\!-multiplier representation commonly available in classical constrained stochastic control, where, under suitable qualification assumptions, the scalar multiplier can be chosen with a definite sign. In our framework, the indeterminacy of its value originates from the recursive coupling introduced by the auxiliary utility system and the resulting time\:\!-inconsistent structure, which break the monotonicity mechanisms that typically ensure sign definiteness. This singular feature not only reflects the intrinsic asymmetry of the recursive formulation, but also highlights the genuinely novel aspect of our contribution, where the constrained and recursive structures are jointly incorporated into a time\:\!-inconsistent equilibrium model.}

\textcolor{black}{Identifying structural conditions that recover sign definiteness---such as monotonicity of the BSDE driver with respect to the utility variable, comparison properties, or convexity/qualification assumptions tailored to the recursive constraint---remains an open and promising direction for further investigation.}

\textcolor{black}{Under suitable assumptions, our analysis also encompasses the case in which the constraint takes the form of a risk constraint. In this context, the additional recursive utility can be naturally interpreted as a dynamic risk measure induced by a $g$\:\!-\:\!expectation, following~\cite{rosazzagianin06}; see also~\cite{gundelweber08},~\cite{imkellerdosreis10},~\cite{yong12}, and~\cite{pennerreveillac15}. More generally, the proposed framework may allow for a broad class of admissible constraints. For instance, one may impose a pointwise constraint on the state process at the terminal time, as in~\cite{jizhou06} and~\cite{zhuo18}. Alternatively, the constraint may involve the expected value of a function of the terminal state, as in~\cite{frankowskazhangzhang19}, where local minimizers, Clarke’s tangent cone, and first- and second-order adjacent cones are employed. For additional insights and related developments, we refer to~\cite{tangli94},~\cite{yong99}, and~\cite{wu13}.}

\textcolor{black}{Research on portfolio optimization problems involving recursive utilities and constraints under time inconsistency is still at a relatively early stage of development. For valuable contributions in this direction, we refer to~\cite{ekelandtaflin05}, which investigates bond portfolios;~\cite{huimkellermuller05} and~\cite{cheriditohu11}, which employ martingale methods;~\cite{pirvu07} and~\cite{morenobrombergpirvureveillac13}, focusing on constant relative risk aversion (CRRA) preferences; and~\cite{horstmorenobromberg08,horstmorenobromberg11}, concerning optimal design.}

\textcolor{black}{To conclude, the final illustrative portfolio problem is amenable to numerical implementation. Its formulation, grounded in the proposed recursive and constrained framework, appears to lend itself to computational treatment through discretization or iterative schemes for the associated adjoint systems. This potential for algorithmic realization underscores the practical significance of the approach and its ability to connect rigorous analysis with realistic decision models. Nonetheless, such aspects are intentionally left beyond the present scope of the study.}

\textcolor{black}{The remainder of the paper is organized as follows. Section~\ref{sec:problems} revisits the notions of admissible controls and subgame\:\!-perfect equilibrium policies, providing the framework for the formulation of the constrained problems. Section~\ref{sec:constrained} presents the main theoretical results, where necessary conditions for the constrained problem are derived using a generalized maximum principle, with detailed discussions of the methodology, including adjoint equations and penalty functionals; notably, these remain valid regardless of the boundedness of the control domain. Section~\ref{sec:application} illustrates the applicability of the theory by means of a portfolio management example, thereby bridging the abstract analysis with a concrete financial setting. Finally, Section~\ref{sec:conclusions} summarizes the main contributions and outlines directions for future research.}


\section{Control problems under constraints}
\label{sec:problems}

\textcolor{black}{To introduce the constrained control problem---Problem~\ref{prob:C}---and set the stage for the main result---Theorem~\ref{th:multipliers}---we briefly recall some key elements from our earlier work~\cite{mastrogiacomotarsia23}, which develops the analysis in the unconstrained setting. This serves both to ensure self\:\!-\:\!containment and to highlight the methodological foundations upon which the present constrained framework is built.}


\textcolor{black}{We consider $T\in\mathopen{]}0,\infty\mathclose{[}$ as a deterministic {\it horizon} and $\rounds{\Omega,\FF,\probP}$ as a complete probability space carrying a one\:\!-dimensional Brownian motion $W = (W(t))_{t\in\mathopen{[}0,T\mathclose{]}}$. Let $\F = (\FF_t)_{t\in\mathopen{[}0,T\mathclose{]}}$ denote the completed filtration generated by $W$\:\!\!, with the assumption that $\FF_T = \FF$\:\!\!, so that the filtered space $\rounds{\Omega,\FF,\F,\probP}$ satisfies the usual conditions. For clarity of exposition, we restrict attention to the one\:\!-dimensional state setting. The multidimensional extension entails no essential difficulty---only notational modifications---and can be handled by adapting the same arguments; see, for instance,~\cite{hu17}.}

\textcolor{black}{We also fix $n \in \N^*$\:\!\! and specify a {\it control domain} $U \in \BB(\R^n) \setminus \Set{\!\emptyset\!}$, that is, a nonempty Borel subset of $\R^n$\:\!\!, and designate an {\it initial instant} $t\in\mathopen{[}0,T\mathclose{[}$. For each $t\in\mathopen{[}0,T\mathclose{[}$, any stochastic differential equation (SDE)---either forward (FSDE) or backward (BSDE)---together with its corresponding adapted solutions, is considered over the domain $\mathopen{[}t,T\mathclose{]}\times\Omega$.}

For a suitable $p \in \mathopen{[}2,\infty\mathclose{[}$\:\! (cf.~\cite[Introduction]{hu17} for details), we denote by $\UUU\bracks{t,T}$ the class of processes $\u(\:\!\bm{\cdot}\:\!) \in \LLped{p}{t,T;\R^n}{\F}$ that are $U$-valued, referring to any $\u(\:\!\bm{\cdot}\:\!)$ of $\UUU\bracks{t,T}$ as an {\it admissible control}. Here, $\LLped{p}{t,T;\R^n}{\F}$ stands for the Banach space of $\R^n$-\:\!valued $(\FF_s)_{s\in\mathopen{[}t,T\mathclose{]}}$-\:\!progressively measurable processes $Z(\:\!\bm{\cdot}\:\!)$ on $\mathopen{[}t,T\mathclose{]}\times\Omega$ satisfying
\[
\E\int_t^T\:\!\!\!{\abs{Z(s)}}^{\:\!p}\:\!ds < \infty \:\! ,
\]
where, as is typically done, $\E\:\!\bracks{\,\bm{\cdot}\,}$ represents the expectation operator with respect to $\probP$\:\!\!, acting on a ($\probP$-\:\!integrable) real\:\!-valued $\FF$-\:\!measurable random variable on $\Omega$.

\textcolor{black}{
Let us fix four deterministic maps
\[
b,\sigma \colon \mathopen{[}0,T\mathclose{]} \times \R \times U \to \R \:\! ,
\]
\[
g \colon \mathopen{[}0,T\mathclose{]} \times \R \times U \times \R \times \R \times \mathopen{[}0,T\mathclose{[} \to \R \:\! ,
\]
and
\[
h \colon \R \times \mathopen{[}0,T\mathclose{[} \to \R \:\! ,
\]
such that the following assumption holds. These conditions are primarily meant to guarantee the well\:\!-posedness of all the stochastic equations arising in the framework developed in the sequel; refer to~\cite{hu17}.
}


\begin{assumption}
\label{ass:assumption}
\textcolor{black}{The functions $b,\sigma,g,h$ are continuous with respect to all their respective variables. For any $t\in\mathopen{[}0,T\mathclose{[}$, there exists a constant $K_t\in\mathopen{]}0,\infty\mathclose{[}$ such that, for every mapping $\phi(s,x,\u,y,z;t)$ among $b(s,x,\u)$, $\sigma(s,x,\u)$, $g(s,x,\u,y,z;t)$, and $h(x;t)$, and for all $s\in\mathopen{[}t,T\mathclose{]}$, $\u \in U$\:\!\!, and $x,y,z\in\R$, the linear growth condition
\[
\abs*{\:\!\phi(s,x,\u,y,z;t)} \:\! \leq K_t\mathopen{\big{(}}1+\abs{x}+\abs{\u}+\abs{y}+\abs{z}\mathclose{\big{)}}
\]
holds. Moreover, the functions $b,\sigma,h$ are of (diﬀerentiability) class $C^{\:\!2}$\:\!\! with respect to the variable $x\in\R$. Their first- and second-order partial derivatives, $b_x,b_{xx},\sigma_x,\sigma_{xx}$, are bounded on $\mathopen{[}t,T\mathclose{]} \times \R \times U$ and continuous with respect to $\rounds{x,\u} \in \R \times U$\:\!\!. Similarly, the derivatives $h_x,h_{xx}$ are bounded and continuous on $\R$. In addition, for each fixed $t$, the function $g(\:\!\bm{\cdot}\,;t)$ is of class $C^{\:\!2}$\:\!\ with respect to $\rounds{x,y,z}\in\R^3$\:\!\!; its gradient $\D g(\:\!\bm{\cdot}\,;t)$ and Hessian matrix $\D^2 \:\!\! g(\:\!\bm{\cdot}\,;t)$, both taken with respect to $\rounds{x,y,z}$, are bounded on $\mathopen{[}t,T\mathclose{]} \times \R \times U \times \R \times \R$ and continuous with respect to $\rounds{x,\u,y,z}\in\R \times U \times \R \times \R$.}
\end{assumption} 

\begin{remark}
\textcolor{black}{Assumption~\ref{ass:assumption} is satisfied by linear maps, on their natural domains, such as those arising in the illustrative example of Section~\ref{sec:application}, which concerns a portfolio optimization problem with constraint. By contrast, this assumption may fail to hold in more general settings. For instance, in the one\:\!-dimensional case ($n = 1$), a ``\:\!toy'' diffusion coefficient $\sigma$ of the form
\[
\sigma(s,x,\u) = s + {\abs{\:\!x\:\!}}^{\:\!1/2}\:\!\! + \u \:\! , \quad (s,x,\u) \in \mathopen{[}0,T\mathclose{]} \times \R \times U ,
\]
does not meet the required smoothness and boundedness properties as the state variable approaches zero, that is, for $x = 0$, regardless of $s$ and $\u$.
}
\end{remark}

Next, we fix a {\it state domain} $A \subseteq \R$, assumed to be a nonempty open interval. For any $t\in\mathopen{[}0,T\mathclose{[}$, $a \in A$, and $\u(\:\!\bm{\cdot}\:\!)\in\UUU\bracks{t,T}$, the FBSDE 
\begin{equation}
\label{eq:fbsde}
\begin{sistema}
dx(s) = b(s,x(s),\u(s))\:\!ds + \sigma(s,x(s),\u(s))\:\!dW(s), \\[0.5ex]
dy(s\:\!;t) = - \:\! g(s,x(s),\u(s),y(s\:\!;t),z(s\:\!;t);t)\:\!ds + z(s\:\!;t)\:\!dW(s), \\[0.5ex]
x(t) = a,\quad y(T\:\!;t) = h(x(T);t),
\end{sistema}
\end{equation} 
represents a {\it recursive stochastic control problem}. If $\rounds{x(\:\!\bm{\cdot}\:\!),y(\:\!\bm{\cdot}\,;t),z(\:\!\bm{\cdot}\,;t)}$ is a solution to~\eqref{eq:fbsde} belonging to $\LLped{2}{t,T;\R}{\F} \times \LLped{2}{\Omega;\CC\rounds{\mathopen{[}t,T\mathclose{]};\R}}{\F} \times \LLped{2}{t,T;\R}{\F}$, then the \textcolor{black}{quadruple} $\rounds{\u(\:\!\bm{\cdot}\:\!),x(\:\!\bm{\cdot}\:\!),y(\:\!\bm{\cdot}\,;t),z(\:\!\bm{\cdot}\,;t)}$ is called {\it admissible}{\color{black}, or equivalently {\it feasible}}. Here, $\LLped{2}{\Omega;\CC\rounds{\mathopen{[}t,T\mathclose{]};\R}}{\F}$ denotes the space of real-\:\!valued, $(\FF_s)_{s\in\mathopen{[}t,T\mathclose{]}}$-\:\!adapted processes $Z(\:\!\bm{\cdot}\:\!)$ on $\mathopen{[}t,T\mathclose{]}\times\Omega$ with ($\probP$-\:\!a.s.) continuous paths and such that
\[
\E\:\!\!\sup_{s\in\mathopen{[}t,T\mathclose{]}} \:\!\! Z^{\:\!2}(s)\:\!\! < \infty \:\! .
\]

If Assumption~\ref{ass:assumption} holds, then for any $t\in\mathopen{[}0,T\mathclose{[}$, $a \in A$, and $\u(\:\!\bm{\cdot}\:\!)\in\UUU\bracks{t,T}$, there exists a unique {\color{black}\emph{strong}} solution
\[
\rounds{x(\:\!\bm{\cdot}\:\!),y(\:\!\bm{\cdot}\,;t),z(\:\!\bm{\cdot}\,;t)} \:\! \in \;\! \LLped{2}{\Omega;\CC\rounds{\mathopen{[}t,T\mathclose{]};\R}}{\F}^2\:\!\! \times \LLped{2}{t,T;\R}{\F}
\]
to~\eqref{eq:fbsde}, as established in Proposition~1 of~\cite{mastrogiacomotarsia23}. Furthermore, $y(t\:\!;t)$ is a deterministic constant (whereas $z(t\:\!;t)$ is generally not); see Remark~4.3 in~\cite{debusschefuhrmantessitore07} and Remark~10 in~\cite{mastrogiacomotarsia23}.

In this case, we define the {\it (dis)utility} or {\it cost functional} $J(\:\!\bm{\cdot}\,;t,a)\colon\UUU\bracks{t,T}\to\R$ by
\begin{equation}
\label{eq:utilityfun}
J(\u(\:\!\bm{\cdot}\:\!);t,a) \doteq y(t\:\!;t) = \E\bracks*{\int_t^T\:\!\!\! g(s,x(s),\u(s),y(s\:\!;t),z(s\:\!;t);t)\,ds + h(x(T);t)}.
\end{equation} 

From now on, we assume that, for any $t\in\mathopen{[}0,T\mathclose{[}$, $\u(\:\!\bm{\cdot}\:\!)\in\UUU\bracks{t,T}$, $a \in A$, and $s\in\mathopen{[}t,T\mathclose{]}$ (and $\probP$-\:\!a.s.),
\begin{equation}
\label{eq:assumptionX}
x(s) \equiv x^{\:\!t,a,\u}(s) \in A \:\! .
\end{equation} 
\textcolor{black}{
Here, we do not interpret the state domain $A$ as a strict state constraint. Rather, it is essentially regarded as a natural range of the generic admissible state process $x(\:\!\bm{\cdot}\:\!)$. In particular, $A$ may be viewed as depending on the length of the time horizon $T$\:\!\!. Ultimately, we may always resort to taking
\[
A = \R \:\! .
\]
We emphasize that the present condition~\eqref{eq:assumptionX} does \emph{not} yet constitute the effective constraint underlying our analysis, which will instead emerge later, grounded in the concept of equilibrium controls (Definition~\ref{def:equilibriumpolicy}) and embodied in condition~\eqref{eq:constraint} (Definition~\ref{def:constraint}).
}


\textcolor{black}{
It is important to highlight the distinct role played by the two time indices $t$ and $s$ in our formulation. The parameter $t \in \mathopen{[}0,T\mathclose{[}$ designates the initial time at which the control problem is posed; together with the initial state $a \in A$, it specifies the starting point of the forward\:\!--\:\!backward system. By contrast, the variable $s \in \mathopen{[}t,T\mathclose{]}$ represents the current time along the trajectories of the state and of the auxiliary processes that will later serve as adjoints.
}

\textcolor{black}{
This asymmetry is fundamental in time\:\!-inconsistent problems, and is precisely what gives rise to an intertemporal misalignment: while $s$ evolves forward dynamically, the performance functional $J(\:\!\bm{\cdot}\,;t,a)$ preserves its dependence on the initial time $t$, thereby capturing the recursive nature of the evaluation. Consequently, the preferences encoded at $t$ do not, in general, coincide with those that a decision\:\!-\:\!maker at a later time $s > t$ would adopt if the problem were reassessed from that vantage point.
}

\textcolor{black}{
In view of the above, we adopt the notion of subgame\:\!-perfect equilibrium policies, as formalized in Definition~3.2 of~\cite{ekelandpirvu08} and subsequently reformulated in Definition~7 of~\cite{mastrogiacomotarsia23}, with only notational adjustments to fit our recursive setting.
}

\begin{definition}[\textcolor{black}{Equilibrium policy/control/pair/quadruple}]
\label{def:equilibriumpolicy}

A {\it (subgame\:\!-perfect) equilibrium policy} is a measurable map $\bm{\Phi} \colon \mathopen{[}0,T\mathclose{]} \times A \to U$ such that, for any $t\in\mathopen{[}0,T\mathclose{[}$ and $a \in A$, there exists a unique $A$\:\!-\:\!valued solution
\[
\bar{x}(\:\!\bm{\cdot}\:\!) \eqnot x^{t,a,\bm{\Phi}}(\:\!\bm{\cdot}\:\!) \in \LLped{2}{t,T;\R}{\F}
\]
to the \textcolor{black}{closed\:\!-\:\!loop} FSDE	
\[
\begin{sistema}
dx(s) = b\bigl{(} \;\!\! s,x(s),\bm{\Phi}(s,x(s)) \;\!\! \bigr{)}\:\!ds + \sigma\bigl{(} \;\!\! s,x(s),\bm{\Phi}(s,x(s)) \;\!\! \bigr{)}\:\!dW(s), \\[0.5ex]
x(t) = a.
\end{sistema}
\]
Moreover, setting
\[
\bar{\u}(\:\!\bm{\cdot}\:\!) \doteq \bm{\Phi}(\:\!\bm{\cdot}\,,\bar{x}(\:\!\bm{\cdot}\:\!)) \in \UUU\bracks{t,T} \:\! ,
\]
it is required that, for any $\u(\:\!\bm{\cdot}\:\!) \in \UUU\bracks{t,T}$,
\begin{equation} 
\label{eq:liminf}
\liminf_{\epsilon\downarrow0}\:\!\frac{J(\bar{\u}^\epsilon(\:\!\bm{\cdot}\:\!);t,a) - J(\bar{\u}(\:\!\bm{\cdot}\:\!);t,a)}{\epsilon} \geq 0 \:\! ,
\end{equation} 
where, as $\epsilon\downarrow0$,
\[
\bar{\u}^\epsilon \:\!\! \doteq \bar{\u} + \rounds{ \u - \bar{\u} } \1_{\:\!\! I_t^{\:\!\epsilon}} 
\]
is the spike (or needle) variation of $\bar{\u}(\:\!\bm{\cdot}\:\!)$ with respect to $\u(\:\!\bm{\cdot}\:\!)$, and
\begin{equation}
\label{eq:Eteps}
I_t^{\:\!\epsilon}\:\!\! \eqnot \mathopen{[}t,t+\epsilon\mathclose{]} \:\! .
\end{equation} 
Here, and throughout the paper, $\1_{\:\!\!E}\rounds{\:\!\bm{\cdot}\:\!}$ denotes the indicator function of a measurable set $E$ belonging to an arbitrary $\sigma$-algebra. \textcolor{black}{By construction, it follows that $\bar{\u}^\epsilon(\:\!\bm{\cdot}\:\!)$ also belongs to $\UUU[t,T]$.} In this terminology, we refer to $\bar{\u}(\:\!\bm{\cdot}\:\!)$ as an {\it equilibrium control} (or {\it strategy}); to $\rounds{\bar{\u}(\:\!\bm{\cdot}\:\!),\bar{x}(\:\!\bm{\cdot}\:\!)}$ as an {\it equilibrium pair}; and to $\rounds{\bar{\u}(\:\!\bm{\cdot}\:\!),\bar{x}(\:\!\bm{\cdot}\:\!),\bar{y}(\:\!\bm{\cdot}\,;t),\bar{z}(\:\!\bm{\cdot}\,;t)}$ as an {\it equilibrium quadruple}. 

\end{definition} 


\textcolor{black}{
As widely recognized in the literature, the uniqueness of an equilibrium policy, control, pair, or quadruple, should it exist, can in general not be guaranteed; see, among others,~\cite{landriaultliliyoung18} and the references therein. In the present framework, the recursive constraint narrows the admissible domain and may thus serve as an implicit mechanism of selection, whereas a systematic treatment of equilibrium multiplicity would require a distinct line of investigation and therefore lies beyond the scope of this paper.
}

\begin{definition}[State constraint]
\label{def:constraint}

\textcolor{black}{
Fix two deterministic maps
\begin{equation}
\label{eq:bmgh}
\bm{g} \colon \mathopen{[}0,T\mathclose{]} \times A \times U \times \R \times \R \times \mathopen{[}0,T\mathclose{[} \to \R \:\! , \qquad \bm{h} \colon A \times \mathopen{[}0,T\mathclose{[} \to \R
\end{equation} 
which satisfies conditions analogous to those in Assumption~\ref{ass:assumption}, with $\bm{g}$ and $\bm{h}$ standing as counterparts of $g$ and $h$ (leaving $b$ and $\sigma$ unchanged). For any $t\in\mathopen{[}0,T\mathclose{[}$ and $a \in A$, let $\bm{\Gamma_{t,a}} \subset \R$ be a nonempty closed interval, referred to as the {\it (dis)utility domain}. For each $\u(\:\!\bm{\cdot}\:\!)\in\UUU\bracks{t,T}$, let $x(\:\!\bm{\cdot}\:\!) = x^{\:\!t,a,\u}(\:\!\bm{\cdot}\:\!)$ be the corresponding state process, and
\[
\bm{y}(\:\!\bm{\cdot}\,;t) \in \LLped{2}{\Omega;\CC\rounds{\mathopen{[}t,T\mathclose{]};\R}}{\F}
\]
be the admissible solution to
\begin{equation}
\label{eq:YZconstraint}
\begin{sistema}
d\:\!\bm{y}(s\:\!;t) = - \:\! \bm{g}(s,x(s),\u(s),\bm{y}(s\:\!;t),\bm{z}(s\:\!;t);t)\:\!ds + \bm{z}(s\:\!;t)\:\!dW(s), \\[0.5ex]
\bm{y}(T\:\!;t) = \bm{h}(x(T);t).
\end{sistema}
\end{equation} 
Define the associated (recursive) utility functional $\bm{J}(\:\!\bm{\cdot}\,;t,a)\colon\UUU\bracks{t,T}\to\R$ by
\begin{equation}
\label{eq:Jconstraint}
\bm{J}(\u(\:\!\bm{\cdot}\:\!);t,a) \doteq \bm{y}(t\:\!;t) = \E\bracks*{\int_t^T\:\!\!\! \bm{g}(s,x(s),\u(s),\bm{y}(s\:\!;t),\bm{z}(s\:\!;t);t)\,ds + \bm{h}(x(T);t)}.
\end{equation} 
We say that the pair {\it $\rounds{\u(\:\!\bm{\cdot}\:\!),x(\:\!\bm{\cdot}\:\!)}$ satisfies the state constraint} (or {\it is subject to the state constraint}) if, and only if,
\begin{equation} 
\label{eq:constraint}
\bm{J}(\u(\:\!\bm{\cdot}\:\!);t,a) \in \bm{\Gamma_{t,a}} \:\! ,
\end{equation} 
and we call condition~\eqref{eq:constraint} the \textit{state constraint for $\rounds{\u(\:\!\bm{\cdot}\:\!),x(\:\!\bm{\cdot}\:\!)}$}. 
}

\end{definition} 

\textcolor{black}{
Concerning the well\:\!-posedness of BSDE~\eqref{eq:YZconstraint} and consequently $\bm{J}(\:\!\bm{\cdot}\,;t,a)$ in~\eqref{eq:Jconstraint}, it is to be observed that the domain of the function $\bm{g}$ as specified in~\eqref{eq:bmgh} is fully consistent with assumption~\eqref{eq:assumptionX}, recalling that this guarantees $x(s) = x^{\:\!t,a,\u}(s) \in A$ ($\probP$-\:\!a.s.) for all admissible choices of the parameters involved.}

\textcolor{black}{The constraint condition~\eqref{eq:constraint} is imposed on $\bm{J}(\:\!\bm{\cdot}\,;t,a)$ in~\eqref{eq:Jconstraint} and not directly on any state process $x(\:\!\bm{\cdot}\:\!)$. The mild abuse of terminology in referring to~\eqref{eq:constraint} in this sense is motivated by~\cite[Chap. 3, Sect. 6]{yongzhou99}, since $\bm{J}(\:\!\bm{\cdot}\,;t,a)$ can be regarded as a real\:\!-valued generalized Bolza\:\!-type functional. Naturally, one should be aware that the condition ultimately constrains the control $\u(\:\!\bm{\cdot}\:\!)$.}

\textcolor{black}{The state constraint~\eqref{eq:constraint} applies only to the pair $\rounds{\u(\:\!\bm{\cdot}\:\!),x(\:\!\bm{\cdot}\:\!)}$ and does not involve the whole quadruple $\rounds{\u(\:\!\bm{\cdot}\:\!),x(\:\!\bm{\cdot}\:\!),y(\:\!\bm{\cdot}\,;t),z(\:\!\bm{\cdot}\,;t)}$. In other words, it restricts the admissibility of the control $\u(\:\!\bm{\cdot}\:\!)$ and the corresponding state process $x(\:\!\bm{\cdot}\:\!)$, while the backward components $y(\:\!\bm{\cdot}\,;t)$ and $z(\:\!\bm{\cdot}\,;t)$ are left unconstrained. This highlights an intrinsic imbalance between the forward and backward parts of the system.}

\textcolor{black}{Alongside its abstract formulation, condition~\eqref{eq:constraint} attains practical relevance in applied settings, as illustrated by the following construction.}

\begin{example}
\label{ex:constraint}
\textcolor{black}{For any fixed $t \in \mathopen{[}0,T\mathclose{[}$, $a \in A$, and $\u(\:\!\bm{\cdot}\:\!) \in \UUU\bracks{t,T}$, the BSDE~\eqref{eq:YZconstraint} together with $\bm{J}(\u(\:\!\bm{\cdot}\:\!);t,a)$ admits a natural interpretation in terms of a scalar (static) risk measure in the sense of~\cite{rosazzagianin06}. More precisely, if $\bm{h}(\:\!\bm{\cdot}\,;t)$ coincides with the identity map on $A$, then BSDE~\eqref{eq:YZconstraint} can equivalently be written as
\[
\begin{sistema}
- \:\! d\:\!\bm{y}(s\:\!;t) = \widetilde{\bm{g}}(s,\bm{y}(s\:\!;t),\bm{z}(s\:\!;t);t)\:\!ds - \bm{z}(s\:\!;t)\:\!dW(s), \\[0.5ex]
\bm{y}(T\:\!;t) = x(T),
\end{sistema}
\]
where $x(T)$ still depends on $\u(\:\!\bm{\cdot}\:\!)$. Thus, following the approach of~\cite{rosazzagianin06}, a constraint condition of the form
\[
\bm{J}(\u(\:\!\bm{\cdot}\:\!);t,a) \equiv \bm{y}(t\:\!;t) \leq \bm{\kappa_{\:\!t,a}} , \quad \bm{\kappa_{\:\!t,a}} \:\!\! \in \R
\]
(that is, condition~\eqref{eq:constraint} with $\bm{\Gamma_{t,a}}$ a right\:\!-\:\!closed ray) amounts to a risk constraint induced by a \emph{$g$\:\!-\:\!expectation}, provided that the driver $\widetilde{\bm{g}}$ satisfies suitable structural conditions. A classical instance is given by the entropic risk measure; see, e.g.,~\cite{frittellirosazzagianin04}.}
\end{example} 

\textcolor{black}{
We are now in a position to articulate the central control problem under consideration. The objective is to describe equilibrium pairs that, besides adhering to the recursive dynamics, also comply with condition~\eqref{eq:constraint}.
}

\begin{problem} 
\label{prob:C}

\textcolor{black}{For any $t\in\mathopen{[}0,T\mathclose{[}$ and $a \in A$, determine {\it necessary conditions} for an equilibrium pair $\rounds{\bar{\u}(\:\!\bm{\cdot}\:\!),\bar{x}^{\:\!t,a,\bar{\u}}(\:\!\bm{\cdot}\:\!)}$ satisfying the state constraint~\eqref{eq:constraint}.}

\end{problem} 

\textcolor{black}{
Our examination of Problem~\ref{prob:C}, along with its resolution yielding a \emph{maximum principle} (variant), will form the core of the next section.
}


\section{Main result: a maximum principle}
\label{sec:constrained}

\textcolor{black}{
We address Problem~\ref{prob:C} by extending to the present setting the analytical procedure developed in~\cite[Chap.~3, Sect.~6]{yongzhou99}, under the assumptions introduced in the previous section. This requires revisiting the classical variational arguments and adjoint structures within a broader recursive and constrained formulation, thereby adapting the stochastic maximum principle to the present equilibrium framework. 
}

\textcolor{black}{
The relevance of the result thus obtained lies in providing a Pontryagin-type formulation for constrained equilibrium policies in a time\:\!-inconsistent recursive setting. In contrast with classical approaches, the analysis does not rely on dynamic programming nor on Hamilton-Jacobi-Bellman (HJB) equations, whose applicability is naturally obstructed by time inconsistency. Instead, it combines variational arguments, adjoint BSDEs, and convex-analytic tools to accommodate, within a unified framework, recursive utilities, state constraints, and subgame\:\!-perfect equilibrium conditions.
}


\subsection{Statement of the maximum principle}
\label{subsec:thmmultipliers}

To proceed, we introduce the generalized Hamiltonian function associated with $b$, $\sigma$, and $g$, following the construction on page~18 of~\cite{hu17} and Definitions~12–13 in~\cite{mastrogiacomotarsia23}. For any ${\color{black}{s \in \mathopen{[}0,T\mathclose{]}}}$, ${\color{black}{\u , \bar{\u} \in U}}$\:\!\!, ${\color{black}{x,y,z,\xi,\theta,\Xi \in \R}}$, $t \in \mathopen{[}0,T\mathclose{[}$, and $\bar{a} \in A$, we define 
\begin{multline}
\label{eq:HH}
\HH(s,x,\u,y,z,\xi,\theta,\Xi\:\! ;t,\bar{a},\bar{\u}) \doteq \xi\:\!b(s,x,\u) + \theta\:\!\sigma(s,x,\u) \\[0.75ex]
+ g(s,x,\u,y,z + \xi\:\!\bracks{\sigma(s,x,\u) - \sigma(s,\bar{a},\bar{\u})};t) \\[0.5ex]
+ \textstyle{\frac{1}{2}}\;\! \Xi\;\! {\bracks{\sigma(s,x,\u) - \sigma(s,\bar{a},\bar{\u})}}^2\:\!\!.
\end{multline} 

\textcolor{black}{For any fixed $t \in \mathopen{[}0,T\mathclose{[}$ and $a \in A$, let $\rounds{\bar{\u}(\:\!\bm{\cdot}\:\!),\bar{x}(\:\!\bm{\cdot}\:\!),\bar{y}(\:\!\bm{\cdot}\,;t),\bar{z}(\:\!\bm{\cdot}\,;t)}$ denote an admissible quadruple. Following the notation in~\cite{hu17}, for any map $\phi(s,x,\u)$ representing $b(s,x,\u)$, $\sigma(s,x,\u)$, or their derivatives, and for $s \in \mathopen{[}t,T\mathclose{]}$ (and $\probP$-\:\!a.s.), we set
\[
\phi(s) \eqnot \phi(s,\bar{x}(s),\bar{\u}(s)) , \quad \delta\phi(s) \eqnot \phi(s,\bar{x}(s),\u(s)) - \phi(s) ,
\]
where $\u(\:\!\bm{\cdot}\:\!) \in \UUU\bracks{t,T}$ is arbitrary. Similarly, for any map $\phi(s,x,\u,y,z;t)$ standing for $f(s,x,\u,y,z;t)$ or its derivatives, and for $s \in \mathopen{[}t,T\mathclose{]}$ (and $\probP$-\:\!a.s.), we write
\[
\phi(s\:\!;t) \eqnot \phi(s,\bar{x}(s),\bar{\u}(s),\bar{y}(s\:\!;t),\bar{z}(s\:\!;t);t) .
\]
}

Hence, we introduce the following two adjoint equations:
\begin{equation}
\label{eq:p}
\begin{sistema}
d\:\!\xi(s\:\!;t) = - \:\! f(s,\xi(s\:\!;t),\theta(s\:\!;t);t)\:\!ds + \theta(s\:\!;t)\:\!dW(s), \\[0.5ex]
\xi(T\:\!;t) = h_x(\bar{x}(T);t),
\end{sistema}
\end{equation} 
where
\begin{equation}
\label{eq:f}
f(s,\xi,\theta\:\!;t) \eqnot \bracks[\big]{b_x(s) + g_z(s\:\!;t)\sigma_x(s) + g_y(s\:\!;t)}\xi + \bracks[\big]{\sigma_x(s) + g_z(s\:\!;t)}\theta + g_x(s\:\!;t) \:\! ,
\end{equation} 
and
\begin{equation}
\label{eq:P}
\begin{sistema}
d\:\!\Xi(s\:\!;t) = - \:\! F(s,\Xi(s\:\!;t),\Theta(s\:\!;t);t)\:\!ds + \Theta(s\:\!;t)\:\!dW(s), \\[0.5ex]
\Xi(T\:\!;t) = h_{xx}(\bar{x}(T);t),
\end{sistema}
\end{equation} 
where
\begin{multline}
\label{eq:F}
F(s,\Xi,\Theta\:\! ;t) \eqnot \bracks*{2 \:\! b_x(s) + {\sigma_x(s)}^2\:\!\! + 2 \:\! g_z(s\:\! ;t)\sigma_x(s) + g_y(s\:\! ;t)}\:\! \Xi \\[0.75ex]
+ \bracks[\big]{2 \:\! \sigma_x(s) + g_z(s\:\!;t)}\Theta+\:\! b_{xx}(s)\xi(s\:\! ;t) + \sigma_{xx}(s)\bracks[\big]{g_z(s\:\!;t)\xi(s\:\! ;t) + \theta(s\:\! ;t)} \\[0.75ex]
+ \big{(}1,\xi(s\:\!;t),\sigma_x(s)\xi(s\:\! ;t) + \theta(s\:\! ;t)\big{)} \:\!\! \cdot \D^2 \:\!\! g(s\:\! ;t) \cdot \:\!\! {\big{(}1,\xi(s\:\!;t),\sigma_x(s)\xi(s\:\! ;t) + \theta(s\:\! ;t)\big{)}}^{\:\!\!\T} \:\!\! .
\end{multline} 
\textcolor{black}{Equations~\eqref{eq:p} and~\eqref{eq:P}, both linear BSDEs, give rise to the first- and second-order adjoint equations/processes associated with the admissible quadruple $\rounds{\bar{\u}(\:\!\bm{\cdot}\:\!),\bar{x}(\:\!\bm{\cdot}\:\!),\bar{y}(\:\!\bm{\cdot}\,;t),\bar{z}(\:\!\bm{\cdot}\,;t)}$. The existence, uniqueness, and regularity of these processes follow from Proposition~3 and Remark~20 in~\cite{mastrogiacomotarsia23}.
}


\begin{notation}
For any $t\in\mathopen{[}0,T\mathclose{[}$, $a \in A$, an admissible quadruple $\rounds{\bar{\u}(\:\!\bm{\cdot}\:\!),\bar{x}(\:\!\bm{\cdot}\:\!),\bar{y}(\:\!\bm{\cdot}\,;t),\bar{z}(\:\!\bm{\cdot}\,;t)}$, $\u \in U$\:\!\!, and $s\in\mathopen{[}t,T\mathclose{]}$, we define, $\probP$-\:\!a.s.,
\begin{multline}
\label{not:HH}
\delta\:\!\HH(s\:\!;t,\u) \eqnot \HH(s,\bar{x}(s),\u,\bar{y}(s\:\!;t),\bar{z}(s\:\!;t),\xi(s\:\!;t),\theta(s\:\!;t),\Xi(s\:\!;t);t,\bar{x}(s),\bar{\u}(s)) \\
	- \HH(s,\bar{x}(s),\bar{\u}(s),\bar{y}(s\:\!;t),\bar{z}(s\:\!;t),\xi(s\:\!;t),\theta(s\:\!;t),\Xi(s\:\!;t);t,\bar{x}(s),\bar{\u}(s)) .
\end{multline} 
\end{notation}

We are now in a position to state the main result of this section, formulated with reference to Definition~\ref{def:constraint} and the related annexes. \textcolor{black}{The proof is given in Subsection~\ref{subsec:thmmultipliersproof} below and is developed through several structured steps.}

\begin{theorem}[Maximum principle] 
\label{th:multipliers}
Suppose that $U$ is bounded. Let $\rounds{\bar{\u}(\:\!\bm{\cdot}\:\!),\bar{x}(\:\!\bm{\cdot}\:\!),\bar{y}(\:\!\bm{\cdot}\,;t),\bar{z}(\:\!\bm{\cdot}\,;t)}$ be an equilibrium quadruple such that the pair $\rounds{\bar{\u}(\:\!\bm{\cdot}\:\!),\bar{x}(\:\!\bm{\cdot}\:\!)}$ is subject to the state constraint:
\begin{equation}
\label{eq:sconstraint}
\bm{J}(\bar{\u}(\:\!\bm{\cdot}\:\!);t,a) = \bar{\bm{y}}(t\:\!;t) \in \bm{\Gamma_{t,a}} \:\! .
\end{equation} 
Then, there exist two \emph{multipliers} $\psi \equiv \psi_{\:\!t,a,\bm{\Gamma_{\:\!\!t,a}}},\bm{\psi} \equiv \bm{\psi}_{\:\!t,a,\bm{\Gamma_{\:\!\!t,a}}} \in \mathopen{[}-1,1\mathclose{]}$, satisfying
\begin{equation}
\label{eq:squaremultipliers}
\psi^2\:\!\! + \bm{\psi}^2\:\!\! = 1,
\end{equation} 
such that the following two conditions hold.
\begin{itemize}[leftmargin=*]

\item[\textbf{1.}]\emph{(Transversality condition)}\textbf{.} For any $\bar{v}\in\bm{\Gamma_{t,a}}$,
\begin{equation}
\label{eq:transversality}
\bm{\psi} \bracks[\big]{\bar{v} - \bm{J}(\bar{\u}(\:\!\bm{\cdot}\:\!);t,a)} \leq 0 \:\! .
\end{equation} 

\item[\textbf{2.}]\emph{ \textcolor{black}{(Maximum condition)}\textbf{.}} Let $\delta\:\!\HH(\:\!\bm{\cdot}\,;t,\textbf{-}\;\!)$ and $\delta\:\!\bm{\HH}(\:\!\bm{\cdot}\,;t,\textbf{-}\;\!)$ be defined with respect to the quadruples $\rounds{\bar{\u}(\:\!\bm{\cdot}\:\!),\bar{x}(\:\!\bm{\cdot}\:\!),\bar{y}(\:\!\bm{\cdot}\,;t),\bar{z}(\:\!\bm{\cdot}\,;t)}$ and $\rounds{\bar{\u}(\:\!\bm{\cdot}\:\!),\bar{x}(\:\!\bm{\cdot}\:\!),\bar{\bm{y}}(\:\!\bm{\cdot}\,;t),\bar{\bm{z}}(\:\!\bm{\cdot}\,;t)}$, respectively. Then, for any $\u \in U$ (and $\probP$-\:\!a.s.),
\begin{equation}
\label{eq:multipliers}
\psi\;\!\delta\:\!\HH(t\:\!;t,\u) + \bm{\psi}\;\!\delta\:\!\bm{\HH}(t\:\!;t,\u) \geq 0 \:\! . 
\end{equation} 

\end{itemize}
\end{theorem} 

\textcolor{black}{
In time\:\!-consistent recursive control problems under constraints, maximum principles can lead, under suitable convexity and well\:\!-posedness assumptions, to a full characterization of optimal policies through an associated coupled forward-backward system; see, in particular,~\cite{elkarouipengquenez01}. In the present framework, however, Theorem~\ref{th:multipliers} is formulated as a necessary condition for a fixed equilibrium quadruple. A direct sufficiency result would require additional structural assumptions under which the Hamiltonian inequality~\eqref{eq:multipliers}, together with the transversality condition~\eqref{eq:transversality}, effectively induces an admissible equilibrium policy and the well\:\!-posedness of the associated FBSDEs. These aspects are closely connected to the existence and uniqueness issues discussed in Remark~\ref{rem:existence} at the end of this subsection, which are intrinsic to the problem.
}

\textcolor{black}{
Against this background, Theorem~\ref{th:multipliers} also conveys several structural implications concerning the interpretation of the multipliers $\psi$ and $\bm{\psi}$. Equality~\eqref{eq:squaremultipliers} ensures that they cannot vanish simultaneously. Whenever $\psi \neq 0$ (that is, $\abs{\:\!\bm{\psi}\:\!} < 1$), we may say that the \emph{qualification condition} holds.
}

\textcolor{black}{
The transversality condition~\eqref{eq:transversality} provides a clear geometric insight: for any $\bar{v}\in\bm{\Gamma_{t,a}}$, the quantities $\bm{\psi}$ and $\bar{v} - \bm{J}(\bar{\u}(\:\!\bm{\cdot}\:\!);t,a)$ have opposite signs. Hence, the sign of $\bm{\psi}$ determines the position of $\bar{\bm{y}}(t;t)$ within the utility domain~$\bm{\Gamma_{t,a}}$ (see~\eqref{eq:sconstraint}). Specifically, if $\bm{\psi} \neq 0$ (i.e., $\abs{\:\!\psi\:\!} < 1$), the following two cases occur.
\[
\bar{\bm{y}}(t\:\!;t) \equiv \bm{J}(\bar{\u}(\:\!\bm{\cdot}\:\!);t,a) =
\begin{cases}
\min\:\!\bm{\Gamma_{t,a}}, & \text{if $\bm{\psi} < 0$,} \\
\max\:\!\bm{\Gamma_{t,a}}, & \text{if $\bm{\psi} > 0$.}
\end{cases}
\]
Conversely,
\begin{equation}
\label{eq:remmultipliers}
\bar{\bm{y}}(t\:\!;t) = \min\:\!\bm{\Gamma_{t,a}}  \quad \Longrightarrow \quad \bm{\psi} \leq 0 \:\! , \qquad \bar{\bm{y}}(t\:\!;t) = \max\:\!\bm{\Gamma_{t,a}} \quad \Longrightarrow \quad \bm{\psi} \geq 0 \:\! .
\end{equation} 
In this sense, $\bm{\psi}$ not only serves as a multiplier but also as a directional indicator of the boundary of the utility domain.
}

\textcolor{black}{
When $\abs{\:\!\psi\:\!} = 1$ ($\bm{\psi} = 0$), inequality~\eqref{eq:multipliers} reduces to the necessary condition part of the maximum principle in~\cite{mastrogiacomotarsia23} (Theorem~1). The state constraint becomes inactive, and the problem coincides with its unconstrained counterpart. More precisely, $\psi = - \:\! 1$ yields $\delta\:\!\HH(t\:\!;t,\u) = 0$ ($\probP$-\:\!a.s.) for every $\u \in U$\:\!\!, which is equivalent to the $\liminf$ term in~\eqref{eq:liminf} being equal to zero.
}

\textcolor{black}{
Further information on the sign of $\bm{\psi}$ can be obtained through the analysis of local perturbations of the Hamiltonians. If there exists $\u^\ast \:\!\! \in U$\:\!\!, which may depend on $t\in\mathopen{[}0,T\mathclose{[}$ or $a \in A$, such that either
\[
\delta\:\!\bm{\HH}(t\:\!;t,\u^\ast \:\!\!) < 0
\]
or
\[
\delta\:\!\bm{\HH}(t\:\!;t,\u^\ast \:\!\!) > 0
\]
($\probP$-\:\!a.s.), then
\[
\bm{\psi} \leq - \;\! \psi \;\! \frac{\delta\:\!\HH(t\:\!;t,\u^\ast \:\!\!)}{\delta\:\!\bm{\HH}(t\:\!;t,\u^\ast \:\!\!)}
\]
or
\[
\bm{\psi} \geq - \;\! \psi \;\! \frac{\delta\:\!\HH(t\:\!;t,\u^\ast \:\!\!)}{\delta\:\!\bm{\HH}(t\:\!;t,\u^\ast \:\!\!)}
\] 
($\probP$-\:\!a.s.), respectively. In this case, if in addition $\delta\:\!\HH(t\:\!;t,\u^\ast \:\!\!) = 0$ ($\probP$-\:\!a.s.), then
\[
\delta\:\!\bm{\HH}(t\:\!;t,\u^\ast \:\!\!) < 0 \quad \Longrightarrow \quad \bm{\psi} \leq 0 \:\! , \qquad \delta\:\!\bm{\HH}(t\:\!;t,\u^\ast \:\!\!) > 0 \quad \Longrightarrow \quad \bm{\psi} \geq 0 \:\! .
\]
Combining this with~\eqref{eq:remmultipliers} yields
\[
\begin{sistema}
\delta\:\!\HH(t\:\!;t,\u^\ast \:\!\!) = 0 \\[0.5ex]
\delta\:\!\bm{\HH}(t\:\!;t,\u^\ast \:\!\!) < 0 \\[0.5ex]
\bar{\bm{y}}(t\:\!;t) = \max\:\!\bm{\Gamma_{t,a}}
\end{sistema}
\!\! \qquad \vee \qquad
\begin{sistema}
\delta\:\!\HH(t\:\!;t,\u^\ast \:\!\!) = 0 \\[0.5ex]
\delta\:\!\bm{\HH}(t\:\!;t,\u^\ast \:\!\!) > 0 \\[0.5ex]
\bar{\bm{y}}(t\:\!;t) = \min\:\!\bm{\Gamma_{t,a}}
\end{sistema}
\qquad \Longrightarrow \qquad \bm{\psi} = 0 \:\! .
\]
}

\textcolor{black}{
Theorem~\ref{th:multipliers} provides necessary optimality conditions formulated with respect to a fixed equilibrium quadruple, rather than addressing the uniqueness of equilibrium controls. A refined notion of local uniqueness will be introduced later, under Assumption~\ref{ass:localuniqueness}; its comprehension naturally benefits from a close reading of the forthcoming proof.
}

\textcolor{black}{
Finally, an important extension of Theorem~\ref{th:multipliers} concerns the case of an unbounded control domain. This setting proves particularly relevant in applications where compactness assumptions are not naturally warranted, thereby enlarging the scope of applicability of the framework; see also~\cite{tangli94, wu13}. The result is stated below, where $\LLped{\infty}{t,T;\R^n}{\F}$ represents the space of $\R^n$-\:\!valued $(\FF_s)_{s\in\mathopen{[}t,T\mathclose{]}}$-\:\!progressively measurable processes $Z(\:\!\bm{\cdot}\:\!)$ on $\mathopen{[}t,T\mathclose{]}\times\Omega$ satisfying
\[
\norma{Z(\:\!\bm{\cdot}\:\!)}_\infty \doteq \inf\Set{\! C\in\mathopen{[}0,\infty\mathclose{[} | \textstyle{\sup_{s\in\mathopen{[}t,T\mathclose{]}} \:\! \abs{Z(s)}} \leq C \text{ ($\probP$-\:\!a.s.)} \!} < \infty \:\! .
\]
}


{\color{black}{
\begin{corollary} 
\label{cor:multipliers}
Let $\rounds{\bar{\u}(\:\!\bm{\cdot}\:\!),\bar{x}(\:\!\bm{\cdot}\:\!),\bar{y}(\:\!\bm{\cdot}\,;t),\bar{z}(\:\!\bm{\cdot}\,;t)}$ be an equilibrium quadruple such that the pair $\rounds{\bar{\u}(\:\!\bm{\cdot}\:\!),\bar{x}(\:\!\bm{\cdot}\:\!)}$ is subject to the state constraint, that is, $\bm{J}(\bar{\u}(\:\!\bm{\cdot}\:\!);t,a) = \bar{\bm{y}}(t\:\!;t) \in \bm{\Gamma_{t,a}}$. Assume also that
\[
\bar{\u}(\:\!\bm{\cdot}\:\!) \in \LLped{\infty}{t,T;\R^n}{\F} \:\! .
\]
Then, there exist $\psi,\bm{\psi}\in\mathopen{[}-1,1\mathclose{]}$, with $\psi^2\:\!\! + \bm{\psi}^2\:\!\! = 1$, such that the following two conditions hold.
\begin{itemize}[leftmargin=*]

\item[\textbf{1.}] For any $\bar{v}\in\bm{\Gamma_{t,a}}$, $\bm{\psi}\bracks[\big]{\bar{v} - \bm{J}(\bar{\u}(\:\!\bm{\cdot}\:\!);t,a)} \leq 0$.

\item[\textbf{2.}] Let $\delta\:\!\HH(\:\!\bm{\cdot}\,;t,\textbf{-}\;\!)$ and $\delta\:\!\bm{\HH}(\:\!\bm{\cdot}\,;t,\textbf{-}\;\!)$ correspond to the quadruples $\rounds{\bar{\u}(\:\!\bm{\cdot}\:\!),\bar{x}(\:\!\bm{\cdot}\:\!),\bar{y}(\:\!\bm{\cdot}\,;t),\bar{z}(\:\!\bm{\cdot}\,;t)}$ and $\rounds{\bar{\u}(\:\!\bm{\cdot}\:\!),\bar{x}(\:\!\bm{\cdot}\:\!),\bar{\bm{y}}(\:\!\bm{\cdot}\,;t),\bar{\bm{z}}(\:\!\bm{\cdot}\,;t)}$, respectively. Then, for any $\u \in U$ (and $\probP$-\:\!a.s.),
\[
\psi\;\!\delta\:\!\HH(t\:\!;t,\u) + \bm{\psi}\;\!\delta\:\!\bm{\HH}(t\:\!;t,\u) \geq 0 \:\! .
\]
\end{itemize}
\end{corollary} 
}}

\textcolor{black}{
The following subsection is devoted to the proof of Theorem~\ref{th:multipliers}, whereas the proof of Corollary~\ref{cor:multipliers} is postponed to Subsection~\ref{subsec:cormultipliers}.
}

{\color{black}{
\begin{remark}
\label{rem:existence}
Theorem~\ref{th:multipliers} provides necessary conditions for equilibrium policies, but does not establish their existence or uniqueness.

At a local level, the existence of a measurable selector satisfying the pointwise optimality condition~\textbf{3} of Theorem~1 in~\cite{mastrogiacomotarsia23} may be guaranteed under standard assumptions, such as compactness of the control domain $U$ and continuity of the Hamiltonian $\HH$ with respect to the control variable $\u$. These conditions imply the non-emptiness of the $\argmin$ and allow for the application of measurable selection theorems.

At a global level, the construction of an equilibrium policy $\bm{\Phi}$ requires solving coupled forward-backward systems, for which existence is typically addressed via fixed-point type arguments, becoming considerably delicate due to the time\:\!-inconsistent recursive structure and the influence of state constraints.

In particular, without an independent existence result, the maximum principle may remain purely formal, in the sense that it characterizes equilibrium candidates without guaranteeing that any equilibrium actually exists.

Concerning uniqueness, additional structural conditions would be required to ensure stability of the forward-backward system and a suitable form of monotonicity or contraction in the associated fixed-point formulation. In the absence of such conditions, multiple equilibria may persist, even under constraints.

A systematic analysis of these issues, including sufficient conditions and uniqueness properties, is beyond the scope of the present paper and is left for future research.
\end{remark}
}}


\subsection{Proof of Theorem~\ref{th:multipliers}}
\label{subsec:thmmultipliersproof}

\textcolor{black}{
The proof unfolds in several steps. We first endow the admissible control space with a suitable complete metric structure, under which the utility functionals are continuous. We then reformulate the constrained problem through an appropriate penalization scheme based on the distance from the constraint set. Ekeland's variational principle is applied to the resulting functional in order to produce an approximate minimizer, whose spike perturbations yield the required first-order variational inequalities.
}

\textcolor{black}{
The final part of the argument relies on a dichotomy, depending on whether the penalized functional remains strictly positive or vanishes along a sequence. In the former case, a normalized representation of the penalized increment leads directly to multiplier estimates, while in the latter case a refined analysis is required to handle degeneracy and recover nontrivial limit multipliers. In each case, suitably normalized multiplier sequences are constructed and shown to converge, ultimately leading to the transversality relation~\eqref{eq:transversality} and the maximum condition~\eqref{eq:multipliers}.
}

\textcolor{black}{
In light of the preceding discussion, we now turn to the technical ingredients underlying the proof.} We begin by restating, for completeness, the notion introduced in Lemma~6.4 of~\cite[Chap.~3, Sect.~6]{yongzhou99}.

\begin{definition}[Metric $\dist$ on $\UUU\bracks{t,T}$]
\label{def:dist}
Fix $t\in\mathopen{[}0,T\mathclose{[}$. Define $\dist \colon \mathopen{\big{(}}\:\!\!\UUU\bracks{t,T}\mathclose{\big{)}}^{\:\!\!2} \! \to \mathopen{[}0,T-t\mathclose{[}$ by
\begin{equation}
\label{eq:dist}
\begin{split}
\dist(\u(\:\!\bm{\cdot}\:\!),\hat{\u}(\:\!\bm{\cdot}\:\!)) &\doteq \rounds{\meas \otimes \probP}\bracks[\big]{\Set{\!(s,\omega)\in\mathopen{[}t,T\mathclose{]}\times\Omega \mid \u(s,\omega)\neq\hat{\u}(s,\omega)\!}} \\[0.75ex]
	&= \E\int_t^T \!\! \1_{\Set{\:\!\!\!\u(\:\!\bm{\cdot}\:\!) \:\! \neq \:\! \hat{\u}(\:\!\bm{\cdot}\:\!)\:\!\!\!}}(s,\omega)\:\!ds \:\! ,
\end{split}
\end{equation} 
where $\meas$ denotes the one\:\!-dimensional Lebesgue measure on $\BB(\R)$.
\end{definition} 

\textcolor{black}{
By construction, Definition~\ref{def:dist} relies on the quotient space of $\UUU\bracks{t,T}$ with respect to the equivalence relation $\u(\:\!\bm{\cdot}\:\!) \sim \tilde{\u}(\:\!\bm{\cdot}\:\!)$ whenever $\dist(\u(\:\!\bm{\cdot}\:\!),\tilde{\u}(\:\!\bm{\cdot}\:\!)) = 0$. This identification, weaker---and therefore more restrictive---than indistinguishability, captures the idea that two admissible controls differing only on a set of vanishing product measure are to be regarded as equivalent. Within this setting, the functional~$\dist$ induces a genuine metric on the resulting space of equivalence classes.
}

\textcolor{black}{
When the control domain~$U$ is bounded, this metric space $\rounds{\UUU\bracks{t,T},\dist}$ is complete, as established in the lemma cited above.
}

We now return to the utility functionals $J$ in~\eqref{eq:utilityfun} and $\bm{J}$ in~\eqref{eq:Jconstraint}, and state the following continuity property, whose proof is deferred to Appendix~\ref{sec:appendix}.

\begin{lemma}
\label{lemma:costcontinuity}
Assume that $U$ is bounded, and fix $t\in\mathopen{[}0,T\mathclose{[}$ and $a \in A$. Then, the utility functionals $J(\:\!\bm{\cdot}\,;t,a),\bm{J}(\:\!\bm{\cdot}\,;t,a)\colon\UUU\bracks{t,T}\to\R$ are continuous with respect to the metric~$\dist$.
\end{lemma} 

\textcolor{black}{
According to the proof of Lemma~\ref{lemma:costcontinuity} (see Appendix~\ref{sec:appendix}), the calculations can be carried out by taking $p = 2$. Moreover, it follows that 
\[
\hat{\u}(\:\!\bm{\cdot}\:\!) \xto{\;\!\dist\;\!}{} \u(\:\!\bm{\cdot}\:\!)
\qquad \Longrightarrow \qquad
\begin{sistema}
\hat{y}(t\:\!;t) \xto{\;\!\R\;\!}{} y(t\:\!;t) \\
\hat{x}(T) \xto{\;\!L^{\:\!\!2}\:\!}{} x(T)
\end{sistema}
\]
as shown, in particular, by~\eqref{eq:claim}.
}

We next recall a preliminary result related to the one\:\!-dimensional version of Lemma~6.5 in~\cite{yongzhou99}.

\begin{lemma} 
\label{lemma:distance}
Let $\Gamma \subseteq \R$ be a nonempty closed interval, and let $\d_\Gamma$ denote the distance function to the set $\Gamma$\:\!\!, that is, the map $\d_\Gamma\colon\R\to\mathopen{[}0,\infty\mathclose{[}$ defined, for each $v\in\R$, by
\begin{equation}
\label{eq:distance}
\d_\Gamma(v) \doteq \inf_{\bar{v}\:\!\in\:\!\Gamma}\;\! \abs*{\:\! v - \bar{v} \:\!} \:\! .
\end{equation} 
Then, $\d_\Gamma$ enjoys the following three properties.

\begin{itemize}[leftmargin=*]

\item[\textbf{1.}] The infimum in~\eqref{eq:distance} is attained as a minimum, and the function $\d_\Gamma$ is convex and $1$-Lipschitz continuous on~$\R$. Moreover,
\begin{equation}
	\label{eq:distancefiber}
	\d_\Gamma^{-1}\:\!\!(0) = \Gamma.
	\end{equation} 

\item[\textbf{2.}] For any $v\in\R$, the subgradient of $\d_\Gamma$ at~$v$, that is, the nonempty subset of~$\R$ 
\begin{equation}
	\label{eq:subgradient}
	\partial\:\!\d_\Gamma(v) = \Set{\! \nu\in\R \;\Big{\vert}\; \text{$\forall\,\hat{v}\in\R$, $\d_\Gamma(\hat{v}) \geq \d_\Gamma(v) + \nu\:\!(\hat{v}-v)$} \!} \:\!\! ,
	\end{equation} 
satisfies
\begin{equation}
	\label{eq:distance1}
	v \notin \Gamma \qquad \Longrightarrow \qquad \partial\:\!\d_\Gamma(v) =
		\begin{cases}
		\Set{\!1\!}\:\!\!, & \text{if $v > \max\:\!\Gamma$,} \\
		\Set{\!-\:\!1\!}\:\!\!, & \text{if $v < \min\:\!\Gamma$.}
		\end{cases}
	\end{equation} 

\item[\textbf{3.}] The squared distance function $\d_\Gamma^2$ is continuously differentiable on~$\R$, and its first derivative is given, for every $v\in\R$, by
\begin{equation}
	\label{eq:distancederivative}
	\frac{d}{dv}\mathopen{\big{(}}\d_\Gamma^2(v)\mathclose{\big{)}} = \:\!
		\begin{cases}
		0, & \text{if $v\in\Gamma$,} \\
		2\:\!\d_\Gamma(v), & \text{if $v > \max\:\!\Gamma$,} \\
		- \:\! 2\:\!\d_\Gamma(v), & \text{if $v < \min\:\!\Gamma$.}
		\end{cases}
	\end{equation} 

\end{itemize}

\end{lemma} 

\begin{notation}
\label{not:distance}
In view of~\eqref{eq:distance1}, for any $v \in \R \setminus \Gamma$, we identify the set $\partial\:\!\d_\Gamma(v)$ with its unique element, that is, $1$ or $-\:\!1$ (a sign). Then, combining~\eqref{eq:distancefiber} and~\eqref{eq:distancederivative}, we can write
\begin{equation*}
\frac{d}{dv}\mathopen{\big{(}}\d_\Gamma^2(\:\!\bm{\cdot}\:\!)\mathclose{\big{)}} = 2\:\!\d_\Gamma(\:\!\bm{\cdot}\:\!)\:\!\partial\:\!\d_\Gamma(\:\!\bm{\cdot}\:\!) \:\! .
\end{equation*} 
\end{notation} 

The content of Lemma~\ref{lemma:distance} may be viewed, in its multidimensional formulation, as a particular instance of {\it Clarke's generalized gradient} theory for locally Lipschitz functions on~$\R^d$\:\!\!. For further details, see~\cite{rockafellar70} and Lemma~2.3 in~\cite[Chap.~3, Sect.~2]{yongzhou99}.

\begin{remark}
\label{rem:distance}
In connection with Lemma~\ref{lemma:distance} (and Notation~\ref{not:distance}), we emphasize that, for any $v\in\R$ and any sequence $(v_\epsilon)_\epsilon$ in $\R$ such that $v_\epsilon \:\!\! \to v$ as $\epsilon\downarrow0$, one has
\[
\d_\Gamma^2(v_\epsilon)-\d_\Gamma^2(v) = \mathopen{\big{[}}2\:\!\d_\Gamma(v)\:\!\partial\:\!\d_\Gamma(v) + o_{\epsilon\downarrow0}(1)\mathclose{\big{]}}(v_\epsilon - v) \:\! ,
\]
or, by continuity,
\[
\d_\Gamma^2(v_\epsilon)-\d_\Gamma^2(v) = \mathopen{\big{[}}2\:\!\d_\Gamma(v_\epsilon)\:\!\partial\:\!\d_\Gamma(v_\epsilon) + \tilde{o}_{\epsilon\downarrow0}(1)\mathclose{\big{]}}(v_\epsilon - v) \:\! ,
\]
for suitable infinitesimals.
\end{remark} 

To prove Theorem~\ref{th:multipliers}, we assume that the control domain~$U$ is bounded and, without loss of generality, that
\begin{equation}
\label{eq:J=0}
J(\bar{\u}(\:\!\bm{\cdot}\:\!);t,a) = 0 \:\! .
\end{equation} 
\textcolor{black}{Observe that~\eqref{eq:J=0} does not directly involve either the utility domain~$\bm{\Gamma_{t,a}}$ or the constraint condition~\eqref{eq:constraint}.}

We work on the metric space $\rounds{\UUU\bracks{t,T},\dist}$ introduced in Definition~\ref{def:dist}, which is complete whenever~$U$ is bounded. The utility functionals $J(\:\!\bm{\cdot}\,;t,a)$ and $\bm{J}(\:\!\bm{\cdot}\,;t,a)$ are continuous with respect to the metric~$\dist$ (see Lemma~\ref{lemma:costcontinuity}).

For any fixed $\rho\downarrow0$, we introduce the {\it penalty functional} $J_\rho(\:\!\bm{\cdot}\,;t,a)\colon\UUU\bracks{t,T}\to\mathopen{[}0,\infty\mathclose{[}$, associated with Problem~\ref{prob:C}, and defined as
\begin{equation}
\label{def:penaltyfunctional}
J_\rho(\u(\:\!\bm{\cdot}\:\!);t,a) \doteq {\braces*{\:\!\!{\rounds{J(\u(\:\!\bm{\cdot}\:\!);t,a)+\rho\:\!}}^2\:\!\! + \d_{\bm{\Gamma_{t,a}}}^2\:\!\!(\bm{J}(\u(\:\!\bm{\cdot}\:\!);t,a))\:\!\!}}^{\:\!\!1/2}
\end{equation} 
where $\d_{\bm{\Gamma_{t,a}}}$ denotes the distance function to the set~$\bm{\Gamma_{t,a}}$, as in~\eqref{eq:distance} of Lemma~\ref{lemma:distance}. If $J(\:\!\bm{\cdot}\,;t,a) \geq 0$ over $\UUU\bracks{t,T}$, i.e., if $\bar{\u}(\:\!\bm{\cdot}\:\!)$ is a classical minimizer of $J(\:\!\bm{\cdot}\,;t,a)$ (see~\eqref{eq:J=0}), then $J_\rho(\:\!\bm{\cdot}\,;t,a) \geq \rho > 0$. Otherwise, the validity of this inequality can no longer be ensured.

It is then straightforward to verify that the functional $J_\rho(\:\!\bm{\cdot}\,;t,a)$ satisfies the following three basic properties.

\begin{itemize}[leftmargin=*]

\item[$\bullet$] $J_\rho(\:\!\bm{\cdot}\,;t,a)$ is continuous with respect to the metric~$\dist$.

\item[$\bullet$] If $\u(\:\!\bm{\cdot}\:\!) = \bar{\u}(\:\!\bm{\cdot}\:\!)$, then
\begin{equation*}
J_\rho(\bar{\u}(\:\!\bm{\cdot}\:\!);t,a) = \rho \:\! ,
\end{equation*}
in light of~\eqref{eq:distancefiber} and~\eqref{eq:J=0}. Moreover,
\[
J_\rho(\bar{\u}(\:\!\bm{\cdot}\:\!);t,a) \leq \! \inf_{\u(\:\!\bm{\cdot}\:\!)\in\;\!\UU\:\!\bracks{t,T}}\:\!\! J_\rho(\u(\:\!\bm{\cdot}\:\!);t,a) + \rho \;\! .
\]

\item[$\bullet$] For any $\u(\:\!\bm{\cdot}\:\!)\in\UUU\bracks{t,T}$,
\begin{equation}
\label{eq:Jrho=0}
J_\rho(\u(\:\!\bm{\cdot}\:\!);t,a) = 0
\qquad \Longleftrightarrow \qquad
\begin{sistema}
J(\u(\:\!\bm{\cdot}\:\!);t,a) = - \:\! \rho \\[1ex]
\bm{J}(\u(\:\!\bm{\cdot}\:\!);t,a) \in \bm{\Gamma_{t,a}}
\end{sistema}
\end{equation} 
(as follows again from~\eqref{eq:distancefiber}).

\end{itemize}

Then, by Corollary~6.3 in~\cite{yongzhou99}---essentially \emph{Ekeland's variational principle} (see~\cite{ekeland74})---there exists $\bar{\u}_\rho(\:\!\bm{\cdot}\:\!)\in\UUU\bracks{t,T}$ such that the following three conditions hold.

\begin{itemize}[leftmargin=*]

\item[\textbf{a.}] $J_\rho(\bar{\u}_\rho(\:\!\bm{\cdot}\:\!);t,a) \leq \rho$\:\!.

\item[\textbf{b.}] $\dist(\bar{\u}(\:\!\bm{\cdot}\:\!),\bar{\u}_\rho(\:\!\bm{\cdot}\:\!)) \leq \sqrt{\rho}$\:\!. Consequently, as $\rho\downarrow0$,
\begin{equation}
	\label{eq:barurho}
	\bar{\u}_\rho(\:\!\bm{\cdot}\:\!) \xto{\;\!\dist\;\!}{} \bar{\u}(\:\!\bm{\cdot}\:\!) \:\! .
	\end{equation} 

\item[\textbf{c.}] For any $\bm{v}(\:\!\bm{\cdot}\:\!)\in\UUU\bracks{t,T}$,
\begin{equation}
	\label{eq:ekelandJrho}
	J_\rho(\bm{v}(\:\!\bm{\cdot}\:\!);t,a) - J_\rho(\bar{\u}_\rho(\:\!\bm{\cdot}\:\!);t,a) \:\! \geq -\:\!\sqrt{\rho}\;\!\dist(\bm{v}(\:\!\bm{\cdot}\:\!),\bar{\u}_\rho(\:\!\bm{\cdot}\:\!)) \:\! .
	\end{equation} 

\end{itemize}

We now fix an arbitrary control $\u(\:\!\bm{\cdot}\:\!)\in\UUU\bracks{t,T}$. For $\epsilon\in\mathopen{]}0,T-t\mathclose{[}$, let $I_t^{\:\!\epsilon}$\:\!\! be as in~\eqref{eq:Eteps}, that is, $I_t^{\:\!\epsilon}\:\!\!=\mathopen{[}t,t+\epsilon\mathclose{]}$, and let $\bar{\u}_\rho^\epsilon(\:\!\bm{\cdot}\:\!)$ denote the spike variation of $\bar{\u}_\rho(\:\!\bm{\cdot}\:\!)$ with respect to $\u(\:\!\bm{\cdot}\:\!)$ and $I_t^{\:\!\epsilon}$\:\!\!. Since $\Set{\!\bar{\u}_\rho^\epsilon(\:\!\bm{\cdot}\:\!) \neq \bar{\u}_\rho(\:\!\bm{\cdot}\:\!)\!} \subseteq I_t^{\:\!\epsilon} \times \Omega$, it follows from~\eqref{eq:dist} that
\begin{equation}
\label{eq:distbarurhoeps}
\dist(\bar{\u}_\rho^\epsilon(\:\!\bm{\cdot}\:\!),\bar{\u}_\rho(\:\!\bm{\cdot}\:\!)) \leq \epsilon \:\! .
\end{equation} 
Hence, as $\epsilon\downarrow0$, the inequality~\eqref{eq:distbarurhoeps} implies that
\begin{equation}
\label{eq:barurhoeps}
\bar{\u}_\rho^\epsilon(\:\!\bm{\cdot}\:\!) \xto{\;\!\dist\;\!}{} \bar{\u}_\rho(\:\!\bm{\cdot}\:\!) \:\! .
\end{equation} 
Moreover, applying~\eqref{eq:ekelandJrho} with $\bm{v}(\:\!\bm{\cdot}\:\!) = \bar{\u}_\rho^\epsilon(\:\!\bm{\cdot}\:\!)$, we obtain
\begin{equation}
\label{eq:ekelandJrhobarurhoeps}
-\:\!\sqrt{\rho}\,\epsilon \leq J_\rho(\bar{\u}_\rho^\epsilon(\:\!\bm{\cdot}\:\!);t,a) - J_\rho(\bar{\u}_\rho(\:\!\bm{\cdot}\:\!);t,a) \:\! .
\end{equation} 

\textcolor{black}{
Recalling the definition~\eqref{def:penaltyfunctional} of the penalized functional $J_\rho(\:\!\bm{\cdot}\,;t,a)$, and observing that~\eqref{eq:Jrho=0} may hold, since $J(\:\!\bm{\cdot}\,;t,a)$ is not a priori nonnegative, we split the proof into two distinct cases. This constitutes the main departure from the classical argument in~\cite{yongzhou99}, where the performance functional can be assumed to be positive. In the present constrained setting, $J_\rho(\:\!\bm{\cdot}\,;t,a)$ may vanish even when the constraint is inactive, which calls for a refined application of Ekeland's variational principle to accommodate weaker sign conditions on the multipliers.
}

\bigskip

\noindent\textbf{Case~I.} $J_\rho(\bar{\u}_\rho(\:\!\bm{\cdot}\:\!);t,a) \neq 0$ (that is, $>0$) for $\rho\downarrow0$\:\!\textbf{.} From~\eqref{def:penaltyfunctional}, we can rewrite~\eqref{eq:ekelandJrhobarurhoeps} as
\begin{equation}
\label{eq:caseI}
\begin{split}
-\:\!\sqrt{\rho}\,\epsilon \:\! &\leq J_\rho(\bar{\u}_\rho^\epsilon(\:\!\bm{\cdot}\:\!);t,a) - J_\rho(\bar{\u}_\rho(\:\!\bm{\cdot}\:\!);t,a) \\[1.5ex]
 &= \frac{J_\rho^{\;\!2}(\bar{\u}_\rho^\epsilon(\:\!\bm{\cdot}\:\!);t,a) - J_\rho^{\;\!2}(\bar{\u}_\rho(\:\!\bm{\cdot}\:\!);t,a)}{J_\rho(\bar{\u}_\rho^\epsilon(\:\!\bm{\cdot}\:\!);t,a) + J_\rho(\bar{\u}_\rho(\:\!\bm{\cdot}\:\!);t,a)} \\[1ex]
 &= \frac{{\rounds{J(\bar{\u}_\rho^\epsilon(\:\!\bm{\cdot}\:\!);t,a)+\rho\:\!}}^2\:\!\! - {\rounds{J(\bar{\u}_\rho(\:\!\bm{\cdot}\:\!);t,a)+\rho\:\!}}^2}{J_\rho(\bar{\u}_\rho^\epsilon(\:\!\bm{\cdot}\:\!);t,a) + J_\rho(\bar{\u}_\rho(\:\!\bm{\cdot}\:\!);t,a)} \\[1ex]
 &+ \frac{\d_{\bm{\Gamma_{t,a}}}^2\:\!\!(\bm{J}(\bar{\u}_\rho^\epsilon(\:\!\bm{\cdot}\:\!);t,a)) - \d_{\bm{\Gamma_{t,a}}}^2\:\!\!(\bm{J}(\bar{\u}_\rho(\:\!\bm{\cdot}\:\!);t,a))}{J_\rho(\bar{\u}_\rho^\epsilon(\:\!\bm{\cdot}\:\!);t,a) + J_\rho(\bar{\u}_\rho(\:\!\bm{\cdot}\:\!);t,a)} 
\end{split}
\end{equation} 
(for any $a_1,a_2\in\R$, if $a_1 + a_2 \neq 0$, then $a_1 - a_2 = \rounds{a_1^2 - a_2^2}/\rounds{a_1 + a_2}$).

We now show that there exist $\psi_\rho^{\:\!\epsilon},\bm{\psi}_\rho^{\:\!\epsilon}\in\R$ such that
\begin{equation}
\label{eq:Krhoeps}
\frac{{\rounds{J(\bar{\u}_\rho^\epsilon(\:\!\bm{\cdot}\:\!);t,a)+\rho\:\!}}^2\:\!\! - {\rounds{J(\bar{\u}_\rho(\:\!\bm{\cdot}\:\!);t,a)+\rho\:\!}}^2}{J_\rho(\bar{\u}_\rho^\epsilon(\:\!\bm{\cdot}\:\!);t,a) + J_\rho(\bar{\u}_\rho(\:\!\bm{\cdot}\:\!);t,a)} = \psi_\rho^{\:\!\epsilon}\:\!\bracks*{J(\bar{\u}_\rho^\epsilon(\:\!\bm{\cdot}\:\!);t,a)-J(\bar{\u}_\rho(\:\!\bm{\cdot}\:\!);t,a)} ,
\end{equation} 
\begin{equation}
\label{eq:bmKrhoeps}
\frac{\d_{\bm{\Gamma_{t,a}}}^2\:\!\!(\bm{J}(\bar{\u}_\rho^\epsilon(\:\!\bm{\cdot}\:\!);t,a)) - \d_{\bm{\Gamma_{t,a}}}^2\:\!\!(\bm{J}(\bar{\u}_\rho(\:\!\bm{\cdot}\:\!);t,a))}{J_\rho(\bar{\u}_\rho^\epsilon(\:\!\bm{\cdot}\:\!);t,a) + J_\rho(\bar{\u}_\rho(\:\!\bm{\cdot}\:\!);t,a)} = \bm{\psi}_\rho^{\:\!\epsilon}\:\!\bracks*{\bm{J}(\bar{\u}_\rho^\epsilon(\:\!\bm{\cdot}\:\!);t,a)-\bm{J}(\bar{\u}_\rho(\:\!\bm{\cdot}\:\!);t,a)} ,
\end{equation} 
and that they are, respectively, of the form
\begin{equation}
\label{eq:psirhoeps}
\psi_\rho^{\:\!\epsilon} = \frac{J(\bar{\u}_\rho(\:\!\bm{\cdot}\:\!);t,a)+\rho}{J_\rho(\bar{\u}_\rho(\:\!\bm{\cdot}\:\!);t,a)} + o_{\epsilon\downarrow0}^{\:\!\rho}(1) \:\! ,
\end{equation} 
\begin{equation}
\label{eq:bmpsirhoeps}
\bm{\psi}_\rho^{\:\!\epsilon} = \frac{\d_{\bm{\Gamma_{t,a}}}\:\!\!(\bm{J}(\bar{\u}_\rho(\:\!\bm{\cdot}\:\!);t,a))\:\!\partial\:\!\d_{\bm{\Gamma_{t,a}}}\:\!\!(\bm{J}(\bar{\u}_\rho(\:\!\bm{\cdot}\:\!);t,a))}{J_\rho(\bar{\u}_\rho(\:\!\bm{\cdot}\:\!);t,a)} + \bm{o}_{\epsilon\downarrow0}^{\:\!\rho}(1) \:\! ,
\end{equation} 
for suitable infinitesimals possibly depending on~$\rho$ (see Notation~\ref{not:distance}).

Indeed, by continuity of $J(\:\!\bm{\cdot}\,;t,a)$, $\bm{J}(\:\!\bm{\cdot}\,;t,a)$ and $J_\rho(\:\!\bm{\cdot}\,;t,a)$ with respect to~$\dist$, together with~\eqref{eq:barurhoeps}, we have
\[
J(\bar{\u}_\rho^\epsilon(\:\!\bm{\cdot}\:\!);t,a) = J(\bar{\u}_\rho(\:\!\bm{\cdot}\:\!);t,a) + o^{\:\!\rho,1}_{\epsilon\downarrow0}(1) \:\! , \qquad J_\rho(\bar{\u}_\rho^\epsilon(\:\!\bm{\cdot}\:\!);t,a) = J_\rho(\bar{\u}_\rho(\:\!\bm{\cdot}\:\!);t,a) + o^{\:\!\rho,2}_{\epsilon\downarrow0}(1)
\]
for suitable infinitesimals. Moreover, by Remark~\ref{rem:distance},
\begin{multline*}
\d_{\bm{\Gamma_{t,a}}}^2\:\!\!(\bm{J}(\bar{\u}_\rho^\epsilon(\:\!\bm{\cdot}\:\!);t,a))-\d_{\bm{\Gamma_{t,a}}}^2\:\!\!(\bm{J}(\bar{\u}_\rho(\:\!\bm{\cdot}\:\! );t,a)) \\[0.75ex]
= \mathopen{\big{[}}2\:\!\d_{\bm{\Gamma_{t,a}}}\:\!\!(\bm{J}(\bar{\u}_\rho(\:\!\bm{\cdot}\:\!);t,a))\:\!\partial\:\!\d_{\bm{\Gamma_{t,a}}}\:\!\!(\bm{J}(\bar{\u}_\rho(\:\!\bm{\cdot}\:\!);t,a)) + \bm{o}^{\:\!\rho,1}_{\epsilon\downarrow0}(1)\mathclose{\big{]}} \:\! \cdot \\[0.75ex]
	\cdot \:\!\! \bracks*{\bm{J}(\bar{\u}_\rho^\epsilon(\:\!\bm{\cdot}\:\!);t,a)-\bm{J}(\bar{\u}_\rho(\:\!\bm{\cdot}\:\!);t,a)}
\end{multline*}
(for infinitesimals as above). Hence, concerning~\eqref{eq:Krhoeps},
\begin{multline*}
{\rounds{J(\bar{\u}_\rho^\epsilon(\:\!\bm{\cdot}\:\!);t,a)+\rho\:\!}}^2\:\!\! - {\rounds{J(\bar{\u}_\rho(\:\!\bm{\cdot}\:\!);t,a)+\rho\:\! }}^2 \:\!\! \\[0.75ex]
= \bracks*{J(\bar{\u}_\rho^\epsilon(\:\!\bm{\cdot}\:\!);t,a)+J(\bar{\u}_\rho(\:\!\bm{\cdot}\:\!);t,a)+2\rho}\bracks*{J(\bar{\u}_\rho^\epsilon(\:\!\bm{\cdot}\:\!);t,a)-J(\bar{\u}_\rho(\:\!\bm{\cdot}\:\! );t,a)} ,
\end{multline*}
and, from what has been previously derived,
\begin{equation*}
\begin{split}
\frac{J(\bar{\u}_\rho^\epsilon(\:\!\bm{\cdot}\:\!);t,a)+J(\bar{\u}_\rho(\:\!\bm{\cdot}\:\!);t,a)+2\rho}{J_\rho(\bar{\u}_\rho^\epsilon(\:\!\bm{\cdot}\:\!);t,a) + J_\rho(\bar{\u}_\rho(\:\!\bm{\cdot}\:\!);t,a)} &= \frac{2\:\!J(\bar{\u}_\rho(\:\!\bm{\cdot}\:\!);t,a)+2\rho+o^{\:\!\rho,1}_{\epsilon\downarrow0}(1)}{2\:\!J_\rho(\bar{\u}_\rho(\:\!\bm{\cdot}\:\!);t,a)+o^{\:\!\rho,2}_{\epsilon\downarrow0}(1)} \\[1ex]
 &\equiv \frac{J(\bar{\u}_\rho(\:\!\bm{\cdot}\:\!);t,a)+\rho}{J_\rho(\bar{\u}_\rho(\:\!\bm{\cdot}\:\!);t,a)} + o^{\:\!\rho}_{\epsilon\downarrow0}(1) \:\! ,
\end{split}
\end{equation*}
which yields~\eqref{eq:psirhoeps}. Likewise, for~\eqref{eq:bmKrhoeps},
\begin{multline*}
\frac{2\:\!\d_{\bm{\Gamma_{t,a}}}\:\!\!(\bm{J}(\bar{\u}_\rho(\:\!\bm{\cdot}\:\!);t,a))\:\!\partial\:\!\d_{\bm{\Gamma_{t,a}}}\:\!\!(\bm{J}(\bar{\u}_\rho(\:\!\bm{\cdot}\:\!);t,a)) + \bm{o}^{\:\!\rho,1}_{\epsilon\downarrow0}(1)}{J_\rho(\bar{\u}_\rho^\epsilon(\:\!\bm{\cdot}\:\! );t,a) + J_\rho(\bar{\u}_\rho(\:\!\bm{\cdot}\:\! );t,a)} \\[1ex]
= \frac{2\:\!\d_{\bm{\Gamma_{t,a}}}\:\!\!(\bm{J}(\bar{\u}_\rho(\:\!\bm{\cdot}\:\!);t,a))\:\!\partial\:\!\d_{\bm{\Gamma_{t,a}}}\:\!\!(\bm{J}(\bar{\u}_\rho(\:\!\bm{\cdot}\:\!);t,a)) + \bm{o}^{\:\!\rho,1}_{\epsilon\downarrow0}(1)}{2\:\ J_\rho(\bar{\u}_\rho(\:\!\bm{\cdot}\:\! );t,a)+o^{\:\!\rho,2}_{\epsilon\downarrow0}(1)} \\[1ex]
\equiv \frac{\d_{\bm{\Gamma_{t,a}}}\:\!\!(\bm{J}(\bar{\u}_\rho(\:\!\bm{\cdot}\:\!);t,a))\:\!\partial\:\!\d_{\bm{\Gamma_{t,a}}}\:\!\!(\bm{J}(\bar{\u}_\rho(\:\!\bm{\cdot}\:\!);t,a))}{J_\rho(\bar{\u}_\rho(\:\!\bm{\cdot}\:\! );t,a)} + \bm{o}^{\:\!\rho}_{\epsilon\downarrow0}(1) \:\! ,
\end{multline*}
leading to~\eqref{eq:bmpsirhoeps}.

The numerator of each fraction characterizing $\psi_\rho^{\:\!\epsilon}$ and $\bm{\psi}_\rho^{\:\!\epsilon}$, unless infinitesimal for $\rho\downarrow0$ by~\eqref{eq:barurho}, may take negative values (see~\eqref{eq:psirhoeps} and~\eqref{eq:bmpsirhoeps}). Nevertheless, when squared,
\[
(\psi_\rho^{\:\!\epsilon})^2\:\!\! + (\bm{\psi}_\rho^{\:\!\epsilon})^2\:\!\! = 1 + \bm{o}^{\:\!\rho,2}_{\epsilon\downarrow0}(1) \:\! ,
\]
as follows from~\eqref{def:penaltyfunctional}. In particular, as $\epsilon\downarrow0$, the sequences $\rounds{\psi_\rho^{\:\!\epsilon}}_\epsilon$ and $\rounds{\bm{\psi}_\rho^{\:\!\epsilon}}_\epsilon$ remain in a compact subset of~$\R$ containing $\mathopen{[}-1,1\mathclose{]}$. Thus, there exist $\psi_\rho,\bm{\psi}_\rho\in\R$ satisfying $(\psi_\rho)^2\:\!\! + (\bm{\psi}_\rho)^2\:\!\! = 1$ and two subsequences---still denoted by $\rounds{\psi_\rho^{\:\!\epsilon}}_\epsilon$ and $\rounds{\bm{\psi}_\rho^{\:\!\epsilon}}_\epsilon$---such that, as $\epsilon\downarrow0$,
\begin{equation}
\label{eq:limpsirhoeps}
\begin{sistema}
\psi_\rho^{\:\!\epsilon} \to \psi_\rho \\[0.75ex]
\bm{\psi}_\rho^{\:\!\epsilon} \to \bm{\psi}_\rho \:\! .
\end{sistema}
\end{equation} 
Analogously, there exist $\psi,\bm{\psi}\in\R$ with
\begin{equation}
\label{eq:squaremultipliersproof}
\psi^2\:\!\! + \bm{\psi}^2\:\!\! = 1
\end{equation} 
and two subsequences---still denoted by $(\psi_\rho)_\rho$ and $(\bm{\psi}_\rho)_\rho$---such that, as $\rho\downarrow0$,
\begin{equation}
\label{eq:limpsirho}
\begin{sistema}
\psi_\rho \to \psi \\[0.75ex]
\bm{\psi}_\rho \to \bm{\psi} \:\! .
\end{sistema}
\end{equation} 
In particular, $\psi,\bm{\psi}\in\mathopen{[}-1,1\mathclose{]}$, so that~\eqref{eq:squaremultipliersproof} coincides with~\eqref{eq:squaremultipliers}.

On the other hand, in connection with~\eqref{eq:transversality}, for every $\bar{v}\in\bm{\Gamma_{t,a}}$ we have
\begin{multline*}
\bm{\psi}_\rho^{\:\!\epsilon}\:\!\mathopen{\big{[}}\bar{v} - \bm{J}(\bar{\u}_\rho(\:\!\bm{\cdot}\:\!);t,a)\mathclose{\big{]}} \\
	= \frac{\d_{\bm{\Gamma_{t,a}}}\:\!\!(\bm{J}(\bar{\u}_\rho(\:\!\bm{\cdot}\:\!);t,a))\:\!\partial\:\!\d_{\bm{\Gamma_{t,a}}}\:\!\!(\bm{J}(\bar{\u}_\rho(\:\!\bm{\cdot}\:\!);t,a))\bracks[\big]{\bar{v} - \bm{J}(\bar{\u}_\rho(\:\!\bm{\cdot}\:\!);t,a)}}{J_\rho(\bar{\u}_\rho(\:\!\bm{\cdot}\:\!);t,a)} + \tilde{\bm{o}}^{\:\!\rho}_{\epsilon\downarrow0}(1)
\end{multline*}
(see~\eqref{eq:bmpsirhoeps}). Here, even if $\bm{J}(\bar{\u}_\rho(\:\!\bm{\cdot}\:\!);t,a) \notin \bm{\Gamma_{t,a}}$,
\begin{multline*}
\partial\:\!\d_{\bm{\Gamma_{t,a}}}\:\!\!(\bm{J}(\bar{\u}_\rho(\:\!\bm{\cdot}\:\! );t,a))\bracks[\big]{\bar{v} - \bm{J}(\bar{\u}_\rho(\:\!\bm{\cdot}\:\! );t,a)} \\[1ex]
\leq \d_{\bm{\Gamma_{t,a}}}\:\!\!(\bar{v}) - \d_{\bm{\Gamma_{t,a}}}\:\!\!(\bm{J}(\bar{\u}_\rho(\:\!\bm{\cdot}\:\! );t,a)) = -\,\d_{\bm{\Gamma_{t,a}}}\:\!\!(\bm{J}(\bar{\u}_\rho(\:\!\bm{\cdot}\:\! );t,a)) \leq 0 \:\! ,
\end{multline*}
where we used the subgradient inequality~\eqref{eq:subgradient} for the distance function with $v = \bm{J}(\bar{\u}_\rho(\:\!\bm{\cdot}\:\!);t,a)$. In particular,
\[
\bm{\psi}_\rho^{\:\!\epsilon}\:\!\mathopen{\big{[}}\bar{v} - \bm{J}(\bar{\u}_\rho(\:\!\bm{\cdot}\:\!);t,a)\mathclose{\big{]}} \leq \tilde{\bm{o}}^{\:\!\rho}_{\epsilon\downarrow0}(1) \:\! ,
\]
and~\eqref{eq:transversality} follows by first letting  $\epsilon\downarrow0$ and then $\rho\downarrow0$ (see~\eqref{eq:barurho},~\eqref{eq:limpsirhoeps}, and~\eqref{eq:limpsirho}).

At this stage, continuing from~\eqref{eq:caseI} through~\eqref{eq:Krhoeps} and~\eqref{eq:bmKrhoeps}, we can write
\begin{equation}
\label{eq:caseIkey}
-\:\!\sqrt{\rho}\,\epsilon \leq \psi_\rho^{\:\!\epsilon}\:\!\bracks*{J(\bar{\u}_\rho^\epsilon(\:\!\bm{\cdot}\:\!);t,a)-J(\bar{\u}_\rho(\:\!\bm{\cdot}\:\!);t,a)} + \bm{\psi}_\rho^{\:\!\epsilon}\:\!\bracks*{\bm{J}(\bar{\u}_\rho^\epsilon(\:\!\bm{\cdot}\:\!);t,a)-\bm{J}(\bar{\u}_\rho(\:\!\bm{\cdot}\:\!);t,a)}.
\end{equation} 
Let $\delta\:\!\HH_\rho(\:\!\bm{\cdot}\,;t,\textbf{-}\;\!)$ correspond to the quadruple $\rounds{\bar{\u}_\rho(\:\!\bm{\cdot}\:\!),\bar{x}_\rho(\:\!\bm{\cdot}\:\!),\bar{y}_\rho(\:\!\bm{\cdot}\,;t),\bar{z}_\rho(\:\!\bm{\cdot}\,;t)}$, and let $\delta\:\!\bm{\HH}_\rho(\:\!\bm{\cdot}\,;t,\textbf{-}\;\!)$ correspond to $\rounds{\bar{\u}_\rho(\:\!\bm{\cdot}\:\!),\bar{x}_\rho(\:\!\bm{\cdot}\:\!),\bar{\bm{y}}_\rho(\:\!\bm{\cdot}\,;t),\bar{\bm{z}}_\rho(\:\!\bm{\cdot}\,;t)}$. Then, by Lemma~2 in~\cite{mastrogiacomotarsia23} with $s = t$, we find
\[
\liminf_{\epsilon\downarrow0}\:\!\frac{J(\bar{\u}_\rho^\epsilon(\:\!\bm{\cdot}\:\!);t,a) - J(\bar{\u}_\rho(\:\!\bm{\cdot}\:\!);t,a)}{\epsilon} = \E\bracks[\big]{\delta\:\!\HH_\rho(t\:\!;t,\u(t))} ,
\]
\[
\liminf_{\epsilon\downarrow0}\:\!\frac{\bm{J}(\bar{\u}_\rho^\epsilon(\:\!\bm{\cdot}\:\!);t,a) - \bm{J}(\bar{\u}_\rho(\:\!\bm{\cdot}\:\!);t,a)}{\epsilon} = \E\bracks[\big]{\delta\:\!\bm{\HH}_\rho(t\:\!;t,\u(t))} ,
\]
and these are indeed true limits (a fact that will be crucial here). Therefore, by~\eqref{eq:limpsirhoeps} and~\eqref{eq:caseIkey},
\[
-\:\!\sqrt{\rho} \;\! \leq \:\! \E\bracks[\big]{\psi_\rho\;\!\delta\:\!\HH_\rho(t\:\!;t,\u(t)) + \bm{\psi}_\rho\;\!\delta\:\!\bm{\HH}_\rho(t\:\!;t,\u(t))} ,
\]
and finally, by~\eqref{eq:barurho} and~\eqref{eq:limpsirho},
\begin{equation}
\label{eq:contradiction}
\E\bracks[\big]{\psi\;\!\delta\:\!\HH(t\:\!;t,\u(t)) + \bm{\psi}\;\!\delta\:\!\bm{\HH}(t\:\!;t,\u(t))} \geq 0 \:\! .
\end{equation} 

We now deduce inequality~\eqref{eq:multipliers}, for any $\u \in U$ (and $\probP$-\:\!a.s.), by contradiction. Indeed, suppose that there exist $t^\ast\:\!\!\in\mathopen{[}0,T\mathclose{[}$, $\u^{\:\!\!\ast}\:\!\! \in U$\:\!\!, and $\NN\in\FF$ with $\probP\bracks{\:\!\NN\:\!} > 0$, such that $\psi\;\!\delta\:\!\HH(t^\ast;t^\ast\:\!\!,\u^{\:\!\!\ast}) + \bm{\psi}\;\!\delta\:\!\bm{\HH}(t^\ast;t^\ast\:\!\!,\u^{\:\!\!\ast}) < 0$ on $\NN$\:\!\!. Then, for any $\u(\:\!\bm{\cdot}\:\!)\in\UUU\bracks{t^\ast\:\!\!,T}$ such that, $\probP$-\:\!a.s.,
\[
\u(t^\ast) =
\begin{cases}
\bar{\u}(t^\ast), & \text{on $\Omega\setminus\NN$\:\!\!,} \\
\u^{\:\!\!\ast\:\!\!}, & \text{on $\NN$\:\!\!,}
\end{cases}
\]
---for instance, the trivial one---we have
\[
\psi\;\!\delta\:\!\HH(t^\ast;t^\ast\:\!\!,\u(t^\ast)) + \bm{\psi}\;\!\delta\:\!\bm{\HH}(t^\ast;t^\ast\:\!\!,\u(t^\ast)) = \bracks[\big]{\psi\;\!\delta\:\!\HH(t^\ast;t^\ast\:\!\!,\u^{\:\!\!\ast}) + \bm{\psi}\;\!\delta\:\!\bm{\HH}(t^\ast;t^\ast\:\!\!,\u^{\:\!\!\ast})}\1_\NN \:\! ,
\]
since $\delta\:\!\HH(t^\ast;t^\ast\:\!\!,\bar{\u}(t^\ast)) = \delta\:\!\bm{\HH}(t^\ast;t^\ast\:\!\!,\bar{\u}(t^\ast)) = 0$. Therefore,
\[
\E\bracks[\big]{\psi\;\!\delta\:\!\HH(t^\ast;t^\ast\:\!\!,\u(t^\ast)) + \bm{\psi}\;\!\delta\:\!\bm{\HH}(t^\ast;t^\ast\:\!\!,\u(t^\ast))} \:\!\! < 0 \:\! ,
\]
because $\probP\bracks{\:\!\NN\:\!} > 0$, which contradicts~\eqref{eq:contradiction}.

\bigskip

\noindent\textbf{Case~II.} For infinitely many $\rho\downarrow0$, $J_\rho(\bar{\u}_\rho(\:\!\bm{\cdot}\:\!);t,a) = 0$, that is,
\[
\begin{sistema}
J(\bar{\u}_\rho(\:\!\bm{\cdot}\:\!);t,a) = - \:\! \rho \\[1ex]
\bm{J}(\bar{\u}_\rho(\:\!\bm{\cdot}\:\!);t,a) \in \bm{\Gamma_{t,a}}
\end{sistema}
\]
(see~\eqref{eq:Jrho=0})\textbf{.} The key point, as in Case~I, is to construct two sequences, $\rounds{\psi_\rho^{\:\!\epsilon}}_\epsilon$ and $\rounds{\bm{\psi}_\rho^{\:\!\epsilon}}_\epsilon$, in~$\R$, bounded as $\epsilon\downarrow0$, such that, similarly to~\eqref{eq:caseI} or~\eqref{eq:caseIkey},
\begin{equation}
\label{eq:caseIIkey}
\begin{split}
J_\rho(\bar{\u}_\rho^\epsilon(\:\!\bm{\cdot}\:\!);t,a) &\equiv J_\rho(\bar{\u}_\rho^\epsilon(\:\!\bm{\cdot}\:\!);t,a)-J_\rho(\bar{\u}_\rho(\:\!\bm{\cdot}\:\!);t,a) \\[0.75ex]
 &\leq \psi_\rho^{\:\!\epsilon}\:\!\bracks*{J(\bar{\u}_\rho^\epsilon(\:\!\bm{\cdot}\:\!);t,a)-J(\bar{\u}_\rho(\:\!\bm{\cdot}\:\!);t,a)} + \bm{\psi}_\rho^{\:\!\epsilon}\:\!\bracks*{\bm{J}(\bar{\u}_\rho^\epsilon(\:\!\bm{\cdot}\:\!);t,a)-\bm{J}(\bar{\u}_\rho(\:\!\bm{\cdot}\:\!);t,a)} \\[0.75ex]
 &\equiv \psi_\rho^{\:\!\epsilon}\:\!\bracks*{J(\bar{\u}_\rho^\epsilon(\:\!\bm{\cdot}\:\!);t,a)+\rho\:\!} + \bm{\psi}_\rho^{\:\!\epsilon}\:\!\bracks*{\bm{J}(\bar{\u}_\rho^\epsilon(\:\!\bm{\cdot}\:\!);t,a)-\bm{J}(\bar{\u}_\rho(\:\!\bm{\cdot}\:\!);t,a)} ,
\end{split}
\end{equation} 
for such values of $\rho\downarrow0$. We shall now analyze this situation by distinguishing two further sub\:\!-\:\!cases, depending on whether or not $\bm{J}(\bar{\u}_\rho(\:\!\bm{\cdot}\:\!);t,a) \in \bm{\Gamma_{t,a}}$.

\bigskip

\noindent\textbf{Subcase~I of Case~II.} For infinitely many of these $\rho\downarrow0$, $\bm{J}(\bar{\u}_\rho(\:\!\bm{\cdot}\:\!);t,a) \in \intpart(\bm{\Gamma_{t,a}})$\textbf{.} For the corresponding values of $\rho\downarrow0$ and for $\epsilon\downarrow0$, we obtain from~\eqref{eq:barurhoeps}
\[
\bm{J}(\bar{\u}_\rho^\epsilon(\:\!\bm{\cdot}\:\!);t,a) \in \intpart(\bm{\Gamma_{t,a}}) .
\]
For such parameters, $J_\rho(\bar{\u}_\rho^\epsilon(\:\!\bm{\cdot}\:\!);t,a) = \abs*{J(\bar{\u}_\rho^\epsilon(\:\!\bm{\cdot}\:\!);t,a)+\rho\:\!}$, since $\intpart(\bm{\Gamma_{t,a}}) \subset \bm{\Gamma_{t,a}}$ (see~\eqref{def:penaltyfunctional}). Therefore, to reach~\eqref{eq:squaremultipliers},~\eqref{eq:transversality}, and~\eqref{eq:caseIIkey}, we simply set
\[
\psi_\rho^{\:\!\epsilon} \eqnot \sgn\rounds{J(\bar{\u}_\rho^\epsilon(\:\!\bm{\cdot}\:\!);t,a)+\rho\:\!} , \qquad \bm{\psi}_\rho^{\:\!\epsilon} \eqnot 0
\]
from which, in the limit, $\abs{\:\!\psi\:\!} = 1$ and $\bm{\psi} = 0$ (as, for any $v\in\R$, $\abs{\:\!v\:\!} = \sgn(v)\:\!v$, or equivalently $\sgn(v)\abs{\:\!v\:\!} = v$).

\bigskip

\noindent\textbf{Subcase~II of Case~II.} For infinitely many of those $\rho\downarrow0$, except for at most a a finite number, $\bm{J}(\bar{\u}_\rho(\:\!\bm{\cdot}\:\!);t,a) \in \partial\:\!\bm{\Gamma_{t,a}}$, and, at the same time, for $\epsilon\downarrow0$,
\[
\bm{J}(\bar{\u}_\rho^\epsilon(\:\!\bm{\cdot}\:\!);t,a) \notin \bm{\Gamma_{t,a}}
\]
(otherwise we can proceed as just seen, possibly passing to a subsequence)\textbf{.} Then, for the corresponding values of $\rho\downarrow0$ and $\epsilon\downarrow0$,
\[
J_\rho(\bar{\u}_\rho^\epsilon(\:\!\bm{\cdot}\:\!);t,a) \neq 0
\]
(that is, $>0$). By~\eqref{def:penaltyfunctional}, and equivalently to~\eqref{eq:caseI} (since $J_\rho(\bar{\u}_\rho(\:\!\bm{\cdot}\:\!);t,a) = 0$), we have
\[
J_\rho(\bar{\u}_\rho^\epsilon(\:\!\bm{\cdot}\:\!);t,a) = \frac{{\rounds{J(\bar{\u}_\rho^\epsilon(\:\!\bm{\cdot}\:\!);t,a)+\rho}}^2\:\!\!}{J_\rho(\bar{\u}_\rho^\epsilon(\:\!\bm{\cdot}\:\!);t,a)} + \frac{\d_{\bm{\Gamma_{t,a}}}^2\:\!\!(\bm{J}(\bar{\u}_\rho^\epsilon(\:\!\bm{\cdot}\:\!);t,a))}{J_\rho(\bar{\u}_\rho^\epsilon(\:\!\bm{\cdot}\:\!);t,a)} \:\! .
\]
Once again, aiming to establish~\eqref{eq:caseIIkey}, we look for $\psi_\rho^{\:\!\epsilon},\bm{\psi}_\rho^{\:\!\epsilon}\in\R$ such that, separately,
\begin{equation}
\label{eq:KrhoepsII}
\frac{{\rounds{J(\bar{\u}_\rho^\epsilon(\:\!\bm{\cdot}\:\!);t,a)+\rho}}^2\:\!\!}{J_\rho(\bar{\u}_\rho^\epsilon(\:\!\bm{\cdot}\:\!);t,a)} \leq \psi_\rho^{\:\!\epsilon}\:\!\bracks*{J(\bar{\u}_\rho^\epsilon(\:\!\bm{\cdot}\:\!);t,a)+\rho\:\!},
\end{equation} 
\begin{equation}
\label{eq:bmKrhoepsII}
\frac{\d_{\bm{\Gamma_{t,a}}}^2\:\!\!(\bm{J}(\bar{\u}_\rho^\epsilon(\:\!\bm{\cdot}\:\!);t,a))}{J_\rho(\bar{\u}_\rho^\epsilon(\:\!\bm{\cdot}\:\!);t,a)} \leq \bm{\psi}_\rho^{\:\!\epsilon}\:\!\bracks*{\bm{J}(\bar{\u}_\rho^\epsilon(\:\!\bm{\cdot}\:\!);t,a)-\bm{J}(\bar{\u}_\rho(\:\!\bm{\cdot}\:\!);t,a)} .
\end{equation} 

Considering~\eqref{eq:KrhoepsII}, we simply take
\[
\psi_\rho^{\:\!\epsilon} \eqnot \frac{J(\bar{\u}_\rho^\epsilon(\:\!\bm{\cdot}\:\!);t,a)+\rho}{J_\rho(\bar{\u}_\rho^\epsilon(\:\!\bm{\cdot}\:\!);t,a)}
\]
(which yields equality in~\eqref{eq:KrhoepsII}); whereas, with regard to~\eqref{eq:bmKrhoepsII}, since
\[
\d_{\bm{\Gamma_{t,a}}}\:\!\!(\bm{J}(\bar{\u}_\rho^\epsilon(\:\!\bm{\cdot}\:\!);t,a)) \leq \partial\:\!\d_{\bm{\Gamma_{t,a}}}\:\!\!(\bm{J}(\bar{\u}_\rho^\epsilon(\:\!\bm{\cdot}\:\!);t,a))\bracks*{\bm{J}(\bar{\u}_\rho^\epsilon(\:\!\bm{\cdot}\:\!);t,a)-\bm{J}(\bar{\u}_\rho(\:\!\bm{\cdot}\:\!);t,a)}
\]
(see~\eqref{eq:subgradient}), we may set
\[
\bm{\psi}_\rho^{\:\!\epsilon} \eqnot \frac{\d_{\bm{\Gamma_{t,a}}}\:\!\!(\bm{J}(\bar{\u}_\rho^\epsilon(\:\!\bm{\cdot}\:\!);t,a))\:\!\partial\:\!\d_{\bm{\Gamma_{t,a}}}\:\!\!(\bm{J}(\bar{\u}_\rho^\epsilon(\:\!\bm{\cdot}\:\!);t,a))}{J_\rho(\bar{\u}_\rho^\epsilon(\:\!\bm{\cdot}\:\!);t,a)} \:\! .
\]
Indeed, on the one hand, concerning~\eqref{eq:squaremultipliers}, we have $(\psi_\rho^{\:\!\epsilon})^2\:\!\! + (\bm{\psi}_\rho^{\:\!\epsilon})^2\:\!\! = 1$; on the other hand, with reference to~\eqref{eq:transversality}, we start from the fact that, for any $\bar{v}\in\bm{\Gamma_{t,a}}$,
\[
\bm{\psi}_\rho^{\:\!\epsilon}\:\!\mathopen{\big{[}}\bar{v} - \bm{J}(\bar{\u}_\rho^\epsilon(\:\!\bm{\cdot}\:\!);t,a)\mathclose{\big{]}} = \frac{\d_{\bm{\Gamma_{t,a}}}\:\!\!(\bm{J}(\bar{\u}_\rho^\epsilon(\:\!\bm{\cdot}\:\!);t,a))\:\!\partial\:\!\d_{\bm{\Gamma_{t,a}}}\:\!\!(\bm{J}(\bar{\u}_\rho^\epsilon(\:\!\bm{\cdot}\:\!);t,a))\bracks[\big]{\bar{v} - \bm{J}(\bar{\u}_\rho^\epsilon(\:\!\bm{\cdot}\:\!);t,a)}}{J_\rho(\bar{\u}_\rho^\epsilon(\:\!\bm{\cdot}\:\!);t,a)} \:\! ,
\]
where $\partial\:\!\d_{\bm{\Gamma_{t,a}}}\:\!\!(\bm{J}(\bar{\u}_\rho^\epsilon(\:\!\bm{\cdot}\:\!);t,a))\bracks[\big]{\bar{v} - \bm{J}(\bar{\u}_\rho^\epsilon(\:\!\bm{\cdot}\:\!);t,a)} < 0$. Hence, we conclude in a perfectly analogous manner to Case~I (see also~\eqref{eq:barurho} and~\eqref{eq:barurhoeps}).

\bigskip

\noindent The proof of Theorem~\ref{th:multipliers} is thus established.

\bigskip

\begin{remark}
Considering~\eqref{eq:KrhoepsII} and~\eqref{eq:bmKrhoepsII}, we may equivalently choose
\[
\psi_\rho^\epsilon = \frac{J(\bar{\u}_\rho^\epsilon(\:\!\bm{\cdot}\:\!);t,a)+\rho}{J_\rho(\bar{\u}_\rho^\epsilon(\:\!\bm{\cdot}\:\!);t,a)} + o^{\:\!\rho}_{\epsilon\downarrow0}(1) \:\! ,
\]
\[
\bm{\psi}_\rho^\epsilon = \frac{\d_{\bm{\Gamma_{t,a}}}\:\!\!(\bm{J}(\bar{\u}_\rho^\epsilon(\:\!\bm{\cdot}\:\!);t,a))\partial\:\!\d_{\bm{\Gamma_{t,a}}}\:\!\!(\bm{J}(\bar{\u}_\rho^\epsilon(\:\!\bm{\cdot}\:\!);t,a))}{J_\rho(\bar{\u}_\rho^\epsilon(\:\!\bm{\cdot}\:\!);t,a)} + \bm{o}^{\:\!\rho}_{\epsilon\downarrow0}(1) \:\! ,
\]
for suitable infinitesimals, not identically zero, such that, respectively,
\[
o^{\:\!\rho}_{\epsilon\downarrow0}(1)\bracks*{J(\bar{\u}_\rho^\epsilon(\:\!\bm{\cdot}\:\!);t,a)+\rho\:\!} \geq 0 , \qquad \bm{o}^{\:\!\rho}_{\epsilon\downarrow0}(1)\bracks*{\bm{J}(\bar{\u}_\rho^\epsilon(\:\!\bm{\cdot}\:\!);t,a)-\bm{J}(\bar{\u}_\rho(\:\!\bm{\cdot}\:\!);t,a)} \geq 0 \:\! .
\]
\end{remark}

Under the forthcoming assumption of \emph{local uniqueness} for the equilibrium quadruple 
$\rounds{\bar{\u}(\:\!\bm{\cdot}\:\!),\bar{x}(\:\!\bm{\cdot}\:\!)$, $\bar{y}(\:\!\bm{\cdot}\:\!;t),\bar{z}(\:\!\bm{\cdot}\:\!;t)}$, where $\rounds{\bar{\u}(\:\!\bm{\cdot}\:\!),\bar{x}(\:\!\bm{\cdot}\:\!)}$ satisfies the state constraint, Case~II of the above proof cannot occur, as it would in the case of a classical minimum point $\bar{\u}(\:\!\bm{\cdot}\:\!)$ for $J(\:\!\bm{\cdot}\,;t,a)$; see~\cite[Chap.~3, Sect.~6]{yongzhou99}.

\medskip

\begin{assumption}
\label{ass:localuniqueness}
There exists $\bar{\delta} > 0$ such that, for any $\hat{\u}(\:\!\bm{\cdot}\:\!)=\hat{\u}_{\:\!\bar{\delta}}(\:\!\bm{\cdot}\:\!)\in\UUU\bracks{t,T}$ with $\bm{J}(\hat{\u}(\:\!\bm{\cdot}\:\!);t,a) \in \bm{\Gamma_{t,a}}$ and $\dist(\bar{\u}(\:\!\bm{\cdot}\:\!),\hat{\u}(\:\!\bm{\cdot}\:\!)) < \bar{\delta}$, if $J(\hat{\u}(\:\!\bm{\cdot}\:\!);t,a) \leq o_{\:\!\bar{\delta}\downarrow0}(1)$ (for a suitable infinitesimal), then $\hat{\u}(\:\!\bm{\cdot}\:\!) = \bar{\u}(\:\!\bm{\cdot}\:\!)$.
\end{assumption} 

\medskip

\noindent Indeed, one can take $\rho\in\mathopen{]}\:\!0,\bar{\delta}^{\:\!2}\mathclose{[}$ and $\hat{\u}(\:\!\bm{\cdot}\:\!) = \bar{\u}_\rho(\:\!\bm{\cdot}\:\!)$ (see~\eqref{eq:J=0}).


\subsection{{\color{black}{Proof of Corollary~\ref{cor:multipliers}}}}
\label{subsec:cormultipliers}


\textcolor{black}{
We proceed by means of a truncation procedure on the control domain, reducing the problem to a sequence of bounded settings to which Theorem~\ref{th:multipliers} is applicable.
}

For any $j\in\N$, consider the Borel subset $\bar{U}_j \in \BB(\R^n) \setminus \Set{\!\emptyset\!}$ of $U$, defined as the intersection of $U$ with the closed ball in $\R^n$\:\!\! centered at the origin and of radius $\norma{\bar{\u}(\:\!\bm{\cdot}\:\!)}_\infty \:\!\! + j$; that is,
\[
\bar{U}_j \eqnot \Set{\! \u \in U \; \big{|} \; \abs*{\u} \leq \norma{\bar{\u}(\:\!\bm{\cdot}\:\!)}_\infty \:\!\! + j \!} \:\!\! .
\]
Next, define the (nonempty) subset $\UUU_j\bracks{t,T}$ of $\UUU\bracks{t,T}$ by
\[
\UUU_j\bracks{t,T} \eqnot \Set{\! \u(\:\!\bm{\cdot}\:\!)\in\LLped{p}{t,T;\R^n\:\!\!}{\F} \; \big{|} \; \forall\; s\in\mathopen{[}t,T\mathclose{]} \text{  {\it and} $\probP$-\:\!a.s.}, \u(s) \in \bar{U}_j \!} \:\!\! .
\]

Observe that $\bar{U}_j$ is bounded for every $j\in \mathbb{N}$. Moreover, for any $s\in\mathopen{[}t,T\mathclose{]}$ and $\probP$-\:\!a.s., one has $\bar{\u}(s) \in \bar{U}_0$; hence, $\bar{\u}(\:\!\bm{\cdot}\:\!)\in\UUU_0\bracks{t,T}$. Since $\bar{U}_j \subseteq \bar{U}_{j+1}$ for all $j\in\N$, it follows that $\UUU_j\bracks{t,T} \subseteq \UUU_{j+1}\bracks{t,T}$. Furthermore,
\begin{equation}
\label{eq:Ucormultipliers}
\bigcup_{j=0}^\infty \:\! \bar{U}_j = U , \qquad \bigcup_{j=0}^\infty \:\! \UUU_j\bracks{t,T} = \UUU\bracks{t,T} \:\! .
\end{equation} 

Now fix $j\in\N$. By the above definitions, Theorem~\ref{th:multipliers} applies with $U$ and $\UUU\bracks{t,T}$ replaced by $\bar{U}_j$ and $\UUU_j\bracks{t,T}$, respectively (noting that $\bar{\u}(\:\!\bm{\cdot}\:\!)\in\UUU_j\bracks{t,T}$, and that the mappings $\HH(\:\!\bm{\cdot}\,;t)$ and $\bm{\HH}(\:\!\bm{\cdot}\,;t)$ can be restricted to $\bar{U}_j$ with respect to the variable $\u$). Thus, there exist $\psi_j,\bm{\psi}_j\in\mathopen{[}-1,1\mathclose{]}$, with $\psi_j^2\:\!\! + \bm{\psi}_j^2\:\!\! = 1$, such that the following two conditions hold.
\begin{itemize}[leftmargin=*]

\item[\textbf{1.}] For any $\bar{v}\in\bm{\Gamma_{t,a}}$, $\bm{\psi}_j\bracks[\big]{\bar{v} - \bm{J}(\bar{\u}(\:\!\bm{\cdot}\:\!);t,a)} \leq 0$.

\item[\textbf{2.}] Let $\delta\:\!\HH(\:\!\bm{\cdot}\,;t,\textbf{-}\;\!)$ and $\delta\:\!\bm{\HH}(\:\!\bm{\cdot}\,;t,\textbf{-}\;\!)$ correspond to the quadruples $\rounds{\bar{\u}(\:\!\bm{\cdot}\:\!),\bar{x}(\:\!\bm{\cdot}\:\!),\bar{y}(\:\!\bm{\cdot}\,;t),\bar{z}(\:\!\bm{\cdot}\,;t)}$ and $\rounds{\bar{\u}(\:\!\bm{\cdot}\:\!),\bar{x}(\:\!\bm{\cdot}\:\!),\bar{\bm{y}}(\:\!\bm{\cdot}\,;t),\bar{\bm{z}}(\:\!\bm{\cdot}\,;t)}$, respectively. Then, for any $\u \in U$ (and $\probP$-\:\!a.s.),
\[
\psi_j\;\!\delta\:\!\HH(t\:\!;t,\u) + \bm{\psi}_j\;\!\delta\:\!\bm{\HH}(t\:\!;t,\u) \geq 0 \:\! .
\]

\end{itemize}

\textcolor{black}{
Finally, by compactness of $\mathopen{[}-1,1\mathclose{]}^2\:\!\!$, one may extract convergent subsequences of $(\psi_j)_j$ and $(\bm{\psi}_j)_j$. Passing to the limit and invoking~\eqref{eq:Ucormultipliers}, we obtain the desired conclusion. This completes the proof of Corollary~\ref{cor:multipliers}, and hence the section.
}



\section{Application to a portfolio optimization problem} 
\label{sec:application}

In this section, we investigate a representative constrained portfolio management problem that serves as an application of the theoretical results derived in the previous section. In particular, we consider a state\:\!-constrained optimization problem over a finite horizon with a non-exponential discounting scheme, such as a hyperbolic one.

The analysis focuses on an intriguing case study concerning the signs of the multipliers $\psi$ and $\bm{\psi}$, {\color{black}under utility functions satisfying a sublinear growth condition. This class includes, in particular, standard CRRA preferences, without requiring a specific parametric form}.

The findings highlight a remarkable connection between the multipliers $\psi$, $\bm{\psi}$, and the first-order process rather than the second-order one. Nevertheless, a certain degree of ambiguity remains regarding the signs of $\psi$ and, more notably, $\bm{\psi}$ themselves, which raises nontrivial interpretative questions in the context of the underlying stochastic dynamics.


With respect to the notations and definitions introduced in Section~5 of~\cite{mastrogiacomotarsia23}, we take the state domain to be
\[
A = \mathopen{]}0,\infty\mathclose{[}
\]
and consider $t\in\mathopen{[}0,T\mathclose{[}$, $a \in A$, two positive processes $\beta(\:\!\bm{\cdot}\,;t),\gamma(\:\!\bm{\cdot}\,;t) \in \LLped{\infty}{t,T;\R}{\F}$, and an admissible investment\:\!-\:\!consumption policy
\[
\rounds{\zeta(\:\!\bm{\cdot}\,,x(\:\!\bm{\cdot}\:\!)),c\:\!(\:\!\bm{\cdot}\,,x(\:\!\bm{\cdot}\:\!))} \in \UUU\bracks{t,T} \:\! ,
\]
taking values in $\R \times \mathopen{[}0,\infty\mathclose{[}$. We then adopt~\eqref{eq:fbsde} in the form
\begin{equation}
\label{eq:applicationuncon}
\begin{sistema}
dx(s) = x(s)\bracks*{\mathopen{\big{(}}r + \mu\:\!\zeta(s,x(s)) - c\:\!(s,x(s))\mathclose{\big{)}}\:\!ds + \sigma\:\!\zeta(s,x(s))\:\!dW(s)} , \\[1ex]
dy(s\:\!;t) = -\:\!\hslash(s\:\!;t)\bracks[\big]{-\:\!\upsilon(c\:\!(s,x(s))x(s)) - \beta(s\:\!;t)y(s\:\!;t) - \gamma(s\:\!;t)z(s\:\!;t)}\:\!ds \\[1ex]
	\hphantom{dY(s\:\!;t) = \:\!} + z(s\:\!;t)\:\!dW(s) \:\! , \\[1ex]
x(t) = a \:\! ,\quad y(T\:\!;t) = -\:\!\hat{\hslash}(T;t)\hat{\upsilon}(x(T)) \:\! ,
\end{sistema}
\end{equation} 
where $r,\mu,\sigma \in \mathopen{]}0,\infty\mathclose{[}$, $\upsilon(\:\!\bm{\cdot}\:\!)$ and $\hat{\upsilon}(\:\!\bm{\cdot}\:\!)$ are two utility functions {\color{black}with sublinear growth}, and $\hslash(\:\!\bm{\cdot}\,;t)$, $\hat{\hslash}(\:\!\bm{\cdot}\,;t)$ are two regular discount functions on $\mathopen{[}t,T\mathclose{]}$. The process $y(\:\!\bm{\cdot}\,;t)$ thereby follows the Uzawa\:\!-type specification (see~\cite{duffieepstein92}).

In order to generate a state constraint for $\rounds{\rounds{\zeta(\:\!\bm{\cdot}\,,x(\:\!\bm{\cdot}\:\!)),c\:\!(\:\!\bm{\cdot}\,,x(\:\!\bm{\cdot}\:\!))},x(\:\!\bm{\cdot}\:\!)}$ as in~\eqref{eq:constraint}, we consider, for any $t\in\mathopen{[}0,T\mathclose{[}$, the parameterized recursive utility system
\begin{equation*}
\begin{sistema}
d\:\!\bm{y}(s\:\!;t) = -\:\!\bm{\hslash}(s\:\!;t)\bracks[\big]{-\:\!k\;\!\bm{\upsilon}(c\:\!(s,x(s))x(s)) - \bm{\beta}(s\:\!;t)\bm{y}(s\:\!;t) - \bm{\gamma}(s\:\!;t)\bm{z}(s\:\!;t)}\:\!ds \\[1ex]
	\hphantom{d\:\!\bm{y}(s\:\!;t) = \:\!} + \bm{z}(s\:\!;t)\:\!dW(s) \:\! , \\[1ex]
\bm{y}(T\:\!;t) = - \:\! \bm{\hat{\hslash}}(T;t)\bm{\hat{\upsilon}}(x(T)) ,
\end{sistema}
\end{equation*} 
where $k \in \mathopen{[}0,\infty\mathclose{[}$, and $\bm{\beta}(\:\!\bm{\cdot}\,;t)$, $\bm{\gamma}(\:\!\bm{\cdot}\,;t) \in \LLped{\infty}{t,T;\R}{\F}$ are positive processes. Moreover, $\bm{\upsilon}(\:\!\bm{\cdot}\:\!)$ and $\bm{\hat{\upsilon}}(\:\!\bm{\cdot}\:\!)$ denote two utility functions {\color{black}with sublinear growth}, while $\bm{\hslash}(\:\!\bm{\cdot}\,;t)$ and $\bm{\hat{\hslash}}(\:\!\bm{\cdot}\,;t)$ are two regular discount functions on $\mathopen{[}t,T\mathclose{]}$.

The previously derived Corollary~\ref{cor:multipliers}, which stems from Theorem~\ref{th:multipliers}, can now be reformulated as follows (see also~\eqref{eq:HH},~\eqref{eq:p},~\eqref{eq:f},~\eqref{eq:P},~\eqref{eq:F}, and~\eqref{not:HH}). \textcolor{black}{We emphasize that the constraint set $\bm{\Gamma_{t,a}}$ denotes a general closed subset of the Euclidean space describing the admissible realizations of the auxiliary recursive utility. In the present formulation, its structure is intentionally left unspecified, so as to encompass a wide range of state\:\!-dependent or expectation\:\!-based constraints.}

\begin{corollary}
\label{cor:applicationconstr}

Let $\rounds{\rounds{\bar{\zeta}(\:\!\bm{\cdot}\,,\bar{x}(\:\!\bm{\cdot}\:\!)),\bar{c}\:\!(\:\!\bm{\cdot}\,,\bar{x}(\:\!\bm{\cdot}\:\!))},\bar{x}(\:\!\bm{\cdot}\:\!),\bar{y}(\:\!\bm{\cdot}\,;t),\bar{z}(\:\!\bm{\cdot}\,;t)}$ be an equilibrium quadruple such that
\[
\bm{J}\big{(}\rounds{\bar{\zeta}(\:\!\bm{\cdot}\,,\bar{x}(\:\!\bm{\cdot}\:\!)),\bar{c}\:\!(\:\!\bm{\cdot}\,,\bar{x}(\:\!\bm{\cdot}\:\!))};t,a\big{)} \in \bm{\Gamma_{t,a}} \:\! .
\]
Assume also that
\[
\rounds{\bar{\zeta}(\:\!\bm{\cdot}\,,\bar{x}(\:\!\bm{\cdot}\:\!)),\bar{c}\:\!(\:\!\bm{\cdot}\,,\bar{x}(\:\!\bm{\cdot}\:\!))} \in \LLped{\infty}{t,T;\R^2}{\F} \:\! .
\]
Then, there exist $\psi,\bm{\psi}\in\mathopen{[}-1,1\mathclose{]}$, with $\psi^2\:\!\! + \bm{\psi}^2\:\!\! = 1$, such that the following two conditions hold.

\begin{itemize}[leftmargin=*]

\item[\textbf{1.}] For any $\bar{v}\in\bm{\Gamma_{t,a}}$, $\bm{\psi}\bracks[\big]{\bar{v} - \bm{J}(\rounds{\bar{\zeta}(\:\!\bm{\cdot}\,,\bar{x}(\:\!\bm{\cdot}\:\!)),\bar{c}\:\!(\:\!\bm{\cdot}\,,\bar{x}(\:\!\bm{\cdot}\:\!))};t,a)} \leq 0$.

\item[\textbf{2.}] Let $\rounds{\:\!\xi(\:\!\bm{\cdot}\,;t),\theta(\:\!\bm{\cdot}\,;t)}$ and $\rounds{\Xi(\:\!\bm{\cdot}\,;t),\Theta(\:\!\bm{\cdot}\,;t)}$ be the adjoint processes corresponding to the quadruple
\[
	\rounds{\rounds{\bar{\zeta}(\:\!\bm{\cdot}\,,\bar{x}(\:\!\bm{\cdot}\:\!)),\bar{c}\:\!(\:\!\bm{\cdot}\,,\bar{x}(\:\!\bm{\cdot}\:\!))},\bar{x}(\:\!\bm{\cdot}\:\!),\bar{y}(\:\!\bm{\cdot}\,;t),\bar{z}(\:\!\bm{\cdot}\,;t)} \:\! ,
	\]
and let $\rounds{\:\!\bm{\xi}(\:\!\bm{\cdot}\,;t),\bm{\theta}(\:\!\bm{\cdot}\,;t)}$ and $\rounds{\bm{\Xi}(\:\!\bm{\cdot}\,;t),\bm{\Theta}(\:\!\bm{\cdot}\,;t)}$ denote the adjoint processes corresponding to the quadruple
\[
	\rounds{\rounds{\bar{\zeta}(\:\!\bm{\cdot}\,,\bar{x}(\:\!\bm{\cdot}\:\!)),\bar{c}\:\!(\:\!\bm{\cdot}\,,\bar{x}(\:\!\bm{\cdot}\:\!))},\bar{x}(\:\!\bm{\cdot}\:\!),\bar{\bm{y}}(\:\!\bm{\cdot}\,;t),\bar{\bm{z}}(\:\!\bm{\cdot}\,;t)} \:\! .
	\]
Then, for any $\rounds{\zeta,c} \in \R \times \mathopen{[}0,\infty\mathclose{[}$ (and $\probP$-\:\!a.s.),
\begin{multline}
	\label{eq:cormultipliers}
	\bracks[\big]{\psi\:\!\xi(t\:\!;t) + \bm{\psi}\:\!\bm{\xi}(t\:\!;t)}\bracks[\big]{\:\!\mu\rounds{\zeta - \bar{\zeta}(t,a)} - \rounds{c - \bar{c}\:\!(t,a)}} + \sigma\bracks[\big]{\psi\:\!\theta(t\:\!;t) + \bm{\psi}\:\!\bm{\theta}(t\:\!;t)}\bracks[\big]{\zeta - \bar{\zeta}(t,a)} \\[0.75ex]
		- \psi\:\!\braces[\bigg]{\:\!\!\frac{1}{a}\bracks[\big]{\upsilon(a\:\!c) - \upsilon(a\:\!\bar{c}\:\!(t,a))} + \sigma\:\!\xi(t\:\!;t)\:\!\gamma(t\:\!;t)\bracks[\big]{\zeta - \bar{\zeta}(t,a)}\:\!\!} \\[0.75ex]
			- \bm{\psi}\:\!\braces[\bigg]{\:\!\!\frac{k}{a}\bracks[\big]{\bm{\upsilon}(a\:\!c) - \bm{\upsilon}(a\:\!\bar{c}\:\!(t,a))} + \sigma\:\!\bm{\xi}(t\:\!;t)\:\!\bm{\gamma}(t\:\!;t)\bracks[\big]{\zeta - \bar{\zeta}(t,a)}\:\!\!} \\[0.75ex]
				+ \frac{a}{2}\:\!\sigma^2\bracks[\big]{\psi\;\!\Xi(t\:\!;t) + \bm{\psi}\;\!\bm{\Xi}(t\:\!;t)}{\bracks[\big]{\zeta - \bar{\zeta}(t,a)}}^2\:\!\! \geq 0 \:\! .
\end{multline} 

\end{itemize}

\end{corollary} 

With respect to Corollary~\ref{cor:applicationconstr}, our focus now shifts to examining the case where
\[
\bm{\psi} \neq 0 \:\! ,
\]
while also allowing for the possibility that $\bm{\psi} = 0$ cannot be excluded.

For this purpose, note that for any $c \in \mathopen{[}0,\infty\mathclose{[}$ (and $\probP$-\:\!a.s.), in~\eqref{eq:cormultipliers}, we encounter a quadratic function with respect to $\zeta - \bar{\zeta}(t,a) \in \R$. This function is nonnegative, implying that it has a nonnegative leading coefficient (i.e., it is convex), namely
\[
\frac{a}{2}\:\!\sigma^2\bracks[\big]{\psi\;\!\Xi(t\:\!;t) + \bm{\psi}\;\!\bm{\Xi}(t\:\!;t)} \:\! ,
\]
and a nonpositive discriminant. We also observe that this function admits a nonzero constant term only if $c \neq \bar{c}\:\!(t,a)$. Therefore, first,
\begin{equation}
\label{eq:disapplicationconstr}
\psi\;\!\Xi(t\:\!;t) + \bm{\psi}\;\!\bm{\Xi}(t\:\!;t) \geq 0 \:\! .
\end{equation} 

Following the spirit of Theorem~3 in~\cite{mastrogiacomotarsia23}, we set
\[
c = \bar{c}\:\!(t,a) \:\! ,
\]
and obtain, for any $k \in \mathopen{[}0,\infty\mathclose{[}$,
\begin{equation}
\label{eq:applicationconstr}
\psi\:\!\braces[\Big]{\! \bracks[\big]{\mu - \sigma\:\!\gamma(t\:\!;t)}\xi(t\:\!;t) + \sigma\:\!\theta(t\:\!;t) \!} + \bm{\psi}\:\!\braces[\Big]{\! \bracks[\big]{\mu - \sigma\:\!\bm{\gamma}(t\:\!;t)}\bm{\xi}(t\:\!;t) + \sigma\:\!\bm{\theta}(t\:\!;t) \!} = 0 \:\! .
\end{equation} 

From here, we separate the analysis into the following cases.

\bigskip

\noindent\textbf{Case~I.} $\bracks[\big]{\mu - \sigma\:\!\bm{\gamma}(t\:\!;t)}\bm{\xi}(t\:\!;t) + \sigma\:\!\bm{\theta}(t\:\!;t) \neq 0$\:\!\textbf{.} From~\eqref{eq:applicationconstr}, it follows that
\[
\bm{\psi} = - \:\! \psi \;\! \frac{\bracks[\big]{\mu - \sigma\:\!\gamma(t\:\!;t)}\xi(t\:\!;t) + \sigma\:\!\theta(t\:\!;t)}{\bracks[\big]{\mu - \sigma\:\!\bm{\gamma}(t\:\!;t)}\bm{\xi}(t\:\!;t) + \sigma\:\!\bm{\theta}(t\:\!;t)} \:\! .
\]
Hence, $\bm{\psi} = 0$ if, and only if, $\bracks[\big]{\mu - \sigma\:\!\gamma(t\:\!;t)}\xi(t\:\!;t) + \sigma\:\!\theta(t\:\!;t) = 0$, which corresponds to~(38) of Theorem~3 in~\cite{mastrogiacomotarsia23} (with $s = t$). Otherwise, since
\[
\psi^2\:\!\! + \bm{\psi}^2\:\!\! = 1 \:\! ,
\]
both $\bm{\psi}$ and $\psi$ are nonzero (that is, $\abs{\:\!\psi\:\!} < 1$ and $\abs{\:\!\bm{\psi}\:\!} < 1$). In particular, using again~\eqref{eq:applicationconstr}, we have
\[
\abs{\:\!\bm{\psi}\:\!} = \frac{\abs*{\bracks[\big]{\mu - \sigma\:\!\gamma(t\:\!;t)}\xi(t\:\!;t) + \sigma\:\!\theta(t\:\!;t)}}{\braces[\Big]{\! \bracks[\big]{\mu - \sigma\:\!\gamma(t\:\!;t)}\xi(t\:\!;t) + \sigma\:\!\theta(t\:\!;t) \!}^{\!2}\:\!\! + \braces[\Big]{\! \bracks[\big]{\mu - \sigma\:\!\bm{\gamma}(t\:\!;t)}\bm{\xi}(t\:\!;t) + \sigma\:\!\bm{\theta}(t\:\!;t) \!}^{\!2}\:\!\!} \;\! ,
\]
\[
\abs{\:\!\psi\:\!} = \frac{\abs*{\bracks[\big]{\mu - \sigma\:\!\bm{\gamma}(t\:\!;t)}\bm{\xi}(t\:\!;t) + \sigma\:\!\bm{\theta}(t\:\!;t)}}{\braces[\Big]{\! \bracks[\big]{\mu - \sigma\:\!\gamma(t\:\!;t)}\xi(t\:\!;t) + \sigma\:\!\theta(t\:\!;t) \!}^{\!2}\:\!\! + \braces[\Big]{\! \bracks[\big]{\mu - \sigma\:\!\bm{\gamma}(t\:\!;t)}\bm{\xi}(t\:\!;t) + \sigma\:\!\bm{\theta}(t\:\!;t) \!}^{\!2}\:\!\!} \;\! .
\]

\bigskip

\noindent\textbf{Case~II.} $\bracks[\big]{\mu - \sigma\:\!\bm{\gamma}(t\:\!;t)}\bm{\xi}(t\:\!;t) + \sigma\:\!\bm{\theta}(t\:\!;t) = 0$\:\!\textbf{.} With reference to the analogues of~\eqref{eq:P} and~\eqref{eq:F} (for the quadruple $\rounds{\rounds{\bar{\zeta}(\:\!\bm{\cdot}\,,\bar{x}(\:\!\bm{\cdot}\:\!)),\bar{c}\:\!(\:\!\bm{\cdot}\,,\bar{x}(\:\!\bm{\cdot}\:\!))},\bar{x}(\:\!\bm{\cdot}\:\!),\bar{\bm{y}}(\:\!\bm{\cdot}\,;t),\bar{\bm{z}}(\:\!\bm{\cdot}\,;t)}$), it is straightforward to verify that $F(\:\!\bm{\cdot}\,,0,0\:\!;t) \geq 0$; since $- \:\! \bm{\hat{\hslash}}(T;t)\bm{\hat{\upsilon}}''(x(T)) > 0$, we obtain
\[
\bm{\Xi}(t\:\!;t) > 0
\]
(see also Remark~17 in~\cite{mastrogiacomotarsia23}). Therefore,~\eqref{eq:disapplicationconstr} becomes
\[
\bm{\psi} \geq - \:\! \psi \;\! \frac{\Xi(t\:\!;t)}{\bm{\Xi}(t\:\!;t)} \:\! .
\]
Hence, if $\bracks[\big]{\mu - \sigma\:\!\gamma(t\:\!;t)}\xi(t\:\!;t) + \sigma\:\!\theta(t\:\!;t) \neq 0$, then, by~\eqref{eq:applicationconstr}, $\psi = 0$ and $\bm{\psi} > 0$ (they cannot both vanish). Otherwise, if $\bracks[\big]{\mu - \sigma\:\!\gamma(t\:\!;t)}\xi(t\:\!;t) + \sigma\:\!\theta(t\:\!;t) = 0$, then, similarly,
\[
\Xi(t\:\!;t) > 0
\]
and, therefore, $\bm{\psi} = 0$ (if and) only if $\psi = 1$ ($\psi = - \:\! 1$ cannot occur); $\psi = 0$ (if and) only if $\bm{\psi} = 1$ ($\bm{\psi} = - \:\! 1$ cannot occur); if $\bm{\psi} < 0$, then $\psi > 0$ (but the converse cannot be affirmed); if $\psi < 0$, then $\bm{\psi} > 0$ (again, the converse cannot be stated). In particular, $\psi$ and $\bm{\psi}$ may both be positive (in which case the above inequality becomes uninformative), but they cannot both be negative.

Consequently, in order to analyze the situations in which it cannot be ruled out that
\[
\bm{\psi} = 0 \:\! ,
\]
we retain the assumptions of Theorem~3 in~\cite{mastrogiacomotarsia23}. Under these assumptions, for any $\rounds{\zeta,c} \in \R \times \mathopen{[}0,\infty\mathclose{[}$ (and $\probP$-\:\!a.s.), inequality~\eqref{eq:cormultipliers} reduces to
\begin{multline*}
	\frac{\psi}{x}\:\!\tilde{\upsilon}(-\:\!\xi(t\:\!;t)) - \psi\braces[\bigg]{\;\!\!\frac{1}{a}\:\!\upsilon(a\:\!c) + c\;\!\xi(t\:\!;t)\:\!\!} \\[0.75ex]
		+ \bm{\psi}\:\!\bm{\xi}(t\:\!;t)\bracks[\big]{\:\!\mu\rounds{\zeta - \bar{\zeta}(t,a)} - \rounds{c - \bar{c}\:\!(t,a)}} + \sigma\:\!\bm{\psi}\:\!\bm{\theta}(t\:\!;t)\bracks[\big]{\zeta - \bar{\zeta}(t,a)} \\[0.75ex]
			- \bm{\psi}\:\!\braces[\bigg]{\:\!\!\frac{k}{a}\bracks[\big]{\bm{\upsilon}(a\:\!c) - \bm{\upsilon}(a\:\!\bar{c}\:\!(t,a))} + \sigma\:\!\bm{\xi}(t\:\!;t)\:\!\bm{\gamma}(t\:\!;t)\bracks[\big]{\zeta - \bar{\zeta}(t,a)}\:\!\!} \\[0.75ex]
				+ \frac{a}{2}\:\!\sigma^2\bracks[\big]{\psi\;\!\Xi(t\:\!;t) + \bm{\psi}\;\!\bm{\Xi}(t\:\!;t)}{\bracks[\big]{\zeta - \bar{\zeta}(t,a)}}^2\:\!\! \geq 0 \:\! ,
\end{multline*}
where $\tilde{\upsilon}(\:\!\bm{\cdot}\:\!)$ denotes the Fenchel\:\!--\:\!Legendre transform of $- \:\! \upsilon(- \;\!\bm{\cdot}\:\!)$. In particular,
\[
\xi(t\:\!;t) < 0 \:\! .
\]

\bigskip

\noindent\textbf{Case~A.} $k = 0$\:\!\textbf{.} By choosing $\zeta = \bar{\zeta}(t,a)$ and $c = 0$, we have
\[
\bm{\xi}(t\:\!;t) < 0
\]
(see~\cite{mastrogiacomotarsia23}). Then, we obtain
\[
\bm{\psi} \leq - \:\! \frac{\psi\;\!\tilde{\upsilon}(-\:\!\xi(t\:\!;t))}{\bm{\xi}(t\:\!;t)\Upsilon(\:\!\xi(t\:\!;t))} \:\! ,
\]
where
\[
\Upsilon \eqnot \big{(} \! -\upsilon' \;\! \big{)}^{-1} .
\]
Therefore, if $\psi \leq 0$, then $\bm{\psi} < 0$ (observing that, if $\psi = 0$, then necessarily $\bm{\psi} \neq 0$). In particular, if $\psi = 0$, then $\bm{\psi} = -\:\!1$. Conversely, if $\bm{\psi} \geq 0$, then $\psi > 0$; in particular, if $\bm{\psi} = 0$, then $\psi = 1$.

More generally, for $\zeta = \bar{\zeta}(t,a)$ and any $c \in \mathopen{[}0,\infty\mathclose{[}$, we have
\[
\bm{\psi} \:\! \sgn\mathopen{\big{(}}\Upsilon(\:\!\xi(t\:\!;t)) - a \:\! c\mathclose{\big{)}} \leq - \:\! \frac{\psi\;\!\braces[\Big]{\tilde{\upsilon}(-\:\!\xi(t\:\!;t)) - \bracks[\big]{\upsilon(a\:\!c) + a \:\! c\;\!\xi(t\:\!;t)}}}{\bm{\xi}(t\:\!;t)\abs[\big]{\:\!\Upsilon(\:\!\xi(t\:\!;t)) - a \:\! c\:\!}} \:\! .
\]
Therefore, for any $c \in \mathopen{]}\bar{c}\:\!(t,a),\infty\mathclose{[}$,
\[
\bm{\psi} \geq - \:\! \frac{\psi\;\!\braces[\Big]{\tilde{\upsilon}(-\:\!\xi(t\:\!;t)) - \bracks[\big]{\upsilon(a\:\!c) + a \:\! c\;\!\xi(t\:\!;t)}\;\!\!}}{\bm{\xi}(t\:\!;t)\bracks[\big]{\:\!\Upsilon(\:\!\xi(t\:\!;t)) - a \:\! c\:\!}} \:\! .
\]
Since $\upsilon(\:\!\bm{\cdot}\:\!)$ is assumed to have sublinear growth, letting $c \uparrow \infty$ yields 
\[
\bm{\psi} \geq - \:\! \frac{\psi \;\! \xi(t\:\!;t)}{\bm{\xi}(t\:\!;t)} \:\! .
\]
It thus appears that, if $\psi \leq 0$, then $\bm{\psi} > 0$; in particular, if $\psi = 0$, then $\bm{\psi} = 1$. Conversely, if $\bm{\psi} \leq 0$, then $\psi > 0$; in particular, if $\bm{\psi} = 0$, then $\psi = 1$. Compared with the previous findings, we can therefore conclude that
\[
\psi > 0
\]
(and, again, that if $\bm{\psi} = 0$, then $\psi = 1$). For $\bm{\psi}$, only the lower bound above is known (that is, $\bm{\psi}$ is greater than or equal to the negative quantity in the inequality). In any case, the qualification condition is satisfied.

\bigskip

\noindent\textbf{Case~B.} $k \neq 0$\:\!\textbf{.} For $\zeta = \bar{\zeta}(t,a)$ and any $c \in \mathopen{[}0,\infty\mathclose{[}$, as before,
\begin{multline*}
\bm{\psi} \:\! \sgn\mathopen{\big{(}}\bm{\xi}(t\:\!;t)\bracks[\big]{\Upsilon(\:\!\xi(t\:\! ;t)) - a \:\! c} + k \:\! \bracks[\big]{\bm{\upsilon}(\Upsilon(\:\!\xi(t\:\! ;t))) - \bm{\upsilon}(a\:\!c)}\mathclose{\big{)}} \\
\geq - \:\! \frac{\psi\;\!\braces[\Big]{\tilde{\upsilon}(-\:\!\xi(t\:\! ;t)) - \bracks[\big]{\upsilon(a\:\!c) + a \:\! c\;\!\xi(t\:\!;t)}}}{\abs*{\:\!\bm{\xi}(t\:\!;t)\bracks[\big]{\Upsilon(\:\!\xi(t\:\! ;t)) - a \:\! c} + k \:\! \bracks[\big]{\bm{\upsilon}(\Upsilon(\:\!\xi(t\:\! ;t))) - \bm{\upsilon}(a\:\!c)}\:\! }} \:\! .
\end{multline*}
Since both $\upsilon(\:\!\bm{\cdot}\:\!)$ and $\bm{\upsilon}(\:\!\bm{\cdot}\:\!)$ exhibit sublinear growth, letting $c \uparrow \infty$ yields 
\[
\bm{\psi} \geq - \:\! \frac{\psi \;\! \xi(t\:\!;t)}{\bm{\xi}(t\:\!;t)} \:\! .
\]
Hence, as before, if $\psi \leq 0$, then $\bm{\psi} > 0$; in particular, if $\psi = 0$, then $\bm{\psi} = 1$. Conversely, if $\bm{\psi} \leq 0$, then $\psi > 0$; in particular, if $\bm{\psi} = 0$, then $\psi = 1$.

However, for $c = 0$,
\[
\bm{\psi} \:\! \sgn\mathopen{\big{(}}\bm{\xi}(t\:\!;t)\Upsilon(\:\!\xi(t\:\!;t)) + k \;\! \bm{\upsilon}(\Upsilon(\:\!\xi(t\:\!;t)))\mathclose{\big{)}} \geq - \:\! \frac{\psi \;\! \tilde{\upsilon}(-\:\!\xi(t\:\!;t))}{\abs[\big]{\:\!\bm{\xi}(t\:\!;t)\Upsilon(\:\!\xi(t\:\!;t)) + k \;\! \bm{\upsilon}(\Upsilon(\:\!\xi(t\:\!;t)))\:\!}} \:\! .
\]

We now distinguish two subcases according to the relative magnitude of the parameter $k$ compared with the other coefficients.

\bigskip

\noindent\textbf{Subcase~I of Case~B.} $k$ is sufficiently large to satisfy
\[
(- \:\! \bm{\xi}(t\:\!;t))\Upsilon(\:\!\xi(t\:\!;t)) < k \;\! \bm{\upsilon}(\Upsilon(\:\!\xi(t\:\!;t))) \:\! \textbf{.}
\]
Then, if $c = 0$,
\[
\bm{\psi} \geq - \:\! \frac{\psi\;\!\tilde{\upsilon}(-\:\!\xi(t\:\!;t))}{\bm{\xi}(t\:\!;t)\Upsilon(\:\!\xi(t\:\!;t)) + k \;\! \bm{\upsilon}(\Upsilon(\:\!\xi(t\:\!;t)))} \:\! .
\]
Therefore, the same conclusions as in the case of an arbitrary $k$ and $c \uparrow \infty$ hold: that is, if $\psi \leq 0$, then $\bm{\psi} > 0$ (in particular, if $\psi = 0$, then $\bm{\psi} = 1$). Conversely, if $\bm{\psi} \leq 0$, then $\psi > 0$ (and in particular, if $\bm{\psi} = 0$, then $\psi = 1$).

\bigskip

\noindent\textbf{Subcase~II of Case~B.} $k$ is sufficiently small to satisfy
\[
(- \:\! \bm{\xi}(t\:\!;t))\Upsilon(\:\!\xi(t\:\!;t)) > k \;\! \bm{\upsilon}(\Upsilon(\:\!\xi(t\:\!;t))) \:\! \textbf{.}
\]
Then, for $c = 0$,
\[
\bm{\psi} \leq - \:\! \frac{\psi\;\!\tilde{\upsilon}(-\:\!\xi(t\:\!;t))}{\bm{\xi}(t\:\!;t)\Upsilon(\:\!\xi(t\:\!;t)) + k \;\! \bm{\upsilon}(\Upsilon(\:\!\xi(t\:\!;t)))}
\]
(consistently with the case $k = 0$ analyzed above). Thus, if $\psi \leq 0$, then $\bm{\psi} < 0$, and if $\psi = 0$, then $\bm{\psi} = -\:\!1$. Conversely, if $\bm{\psi} \geq 0$, then $\psi > 0$; in particular, if $\bm{\psi} = 0$, then $\psi = 1$. Compared with the results obtained above, we can conclude only that
\[
\psi > 0 \:\! .
\]

\textcolor{black}{
While an explicit closed\:\!-\:\!form expression of the equilibrium strategy is not attainable, the formulation may allow for a numerical computation of the constrained equilibrium through standard discretization or backward-iteration schemes applied to the adjoint system. A natural benchmark for such computation is provided by the unconstrained recursive equilibrium derived in~\cite{mastrogiacomotarsia23}, and, when relevant, by the time\:\!-consistent optimal control under exponential discounting. Comparing these outcomes would help quantify the impact of the constraint on consumption and portfolio allocation, thereby reinforcing the methodological relevance of the present approach.
}

\section{Concluding remarks}
\label{sec:conclusions}

\textcolor{black}{This paper has extended the analysis of time\:\!-inconsistent recursive stochastic control problems to the constrained setting. Building on previous contributions, most notably~\cite{ekelandlazrak06, ekelandpirvu08, ekelandmbodjipirvu12, hu17, mastrogiacomotarsia23}, we have derived a Pontryagin-type maximum principle characterizing subgame\:\!-perfect equilibria of closed\:\!-\:\!loop type. The resulting framework introduces a constraint linked to a utility functional defined through an additional recursive utility system. To address this setting, the analysis combines Ekeland's variational principle with a fundamental result on the distance function relative to a closed subset of the Euclidean space.}

\textcolor{black}{The novelty of the contribution is mainly analytical: several existing tools are integrated into a coherent framework that accommodates both recursive utilities and constraints under time inconsistency. Relying on the techniques of~\cite{yongzhou99} and~\cite{hu17}, and crucially on Lemma~2 in~\cite{mastrogiacomotarsia23}, the proof of the main result, Theorem~\ref{th:multipliers}, addresses a number of sub\:\!-\:\!cases not covered by the classical theory.}

\textcolor{black}{This theorem establishes the existence of generalized Lagrange\:\!-\:\!multipliers in a setting where the scalar multiplier~$\psi$ cannot be guaranteed to have a fixed sign, preventing the direct embedding of the multipliers~$\psi$ and~$\bm{\psi}$ into the Hamiltonian expansion. Nevertheless, the formulation proves effective in financial applications, encompassing both trading and consumption models, and offers a rigorous foundation for constrained recursive equilibria.}

\textcolor{black}{In contrast to, for instance,~\cite{hernandezpossamai23, hernandezpossamai24}, which rely on an extended dynamic programming principle and characterize equilibria via an infinite BSDE system analogous to HJB equations, our approach is variational in nature and follows a maximum-principle route. This methodological distinction highlights the flexibility of the present framework, which accommodates recursive utilities used as constraints through penalization and captures their equilibrium structure via adjoint equations.}

\textcolor{black}{These results complement and extend those of~\cite{mastrogiacomotarsia23}, and are specific to the constrained-recursive, time\:\!-inconsistent formulation introduced here, while suggesting several directions for further research, including the analytical characterization of the constrained Hamiltonian and potential extensions to mean-field or multi-agent settings.}

\textcolor{black}{The present framework deliberately relies on regularity assumptions that ensure well\:\!-posedness under globally $C^{\:\!2}$\:\!\! coefficients in the state variable. As a consequence, diffusion dynamics featuring square\:\!-\:\!root or similarly nonlinear behaviour in at least one component---such as those underlying the Cox\:\!--\:\!Ingersoll\:\!--\:\!Ross short\:\!-\:\!rate model (cf.~\cite{coxingersollross85}) and the Heston volatility model (cf.~\cite{heston93})---fall outside the current setting. Nevertheless, relaxing these smoothness requirements offers a promising avenue for further investigation.} 

\bmhead{Acknowledgements}

\textcolor{black}{The authors gratefully acknowledge the support of the GNAMPA (INdAM) group during the development of the research project that resulted in this paper.} 


\section*{Declarations}

The two authors declare that no funds or grants were received during the preparation of this manuscript. They have no relevant financial or non\:\!-\:\!financial interests to disclose. Both authors contributed to the study's conception and design, and all read and approved the final manuscript.

%
%
%

%
%
%
%


\begin{appendices}

\section{Proof of Lemma~\ref{lemma:costcontinuity}}
\label{sec:appendix}

It suffices to establish the result for $J(\:\!\bm{\cdot}\,;t,a)$. Fix $\u(\:\!\bm{\cdot}\:\!),\hat{\u}(\:\!\bm{\cdot}\:\!)\in\UUU\bracks{t,T}$ and consider the corresponding admissible quadruples $\rounds{\u(\:\!\bm{\cdot}\:\!),x(\:\!\bm{\cdot}\:\!),y(\:\!\bm{\cdot}\,;t),z(\:\!\bm{\cdot}\,;t)}$ and $\rounds{\hat{\u}(\:\!\bm{\cdot}\:\!),\hat{x}(\:\!\bm{\cdot}\:\!),\hat{y}(\:\!\bm{\cdot}\,;t),\hat{z}(\:\!\bm{\cdot}\,;t)}$. By Lemma~1 in~\cite{mastrogiacomotarsia23} and Jensen's inequality (in both its continuous and discrete forms),
\begin{equation*}
\begin{split}
\hspace{0.001ex} & {\abs[\big]{J(\u(\:\!\bm{\cdot}\:\!);t,a) - J(\hat{\u}(\:\!\bm{\cdot}\:\!);t,a)}}^{\:\!p}\:\!\! \\[3ex]
 &= {\abs*{y(t\:\!;t) - \hat{y}(t\:\!;t)}}^{\:\!p}\:\!\! \\[3ex]
 &\lesssim \E{\abs*{h(x(T);t) - h(\hat{x}(T);t)}}^{\:\!p}\:\!\! \\[1.5ex]
 &+ \E\:\!\!{\mathopen{\bigg{(}}\int_t^T\!\abs*{g(s,x(s),\u(s),y(s\:\!;t),z(s\:\!;t);t) - g(s,\hat{x}(s),\hat{\u}(s),y(s\:\!;t),z(s\:\!;t);t)}\:\!ds\mathclose{\bigg{)}}}^{\!p\:\!}\:\!\! \\[0.5ex]
 &\lesssim \E{\abs*{x(T) - \hat{x}(T)}}^{\:\!p}\:\!\! + \E\:\!\!{\mathopen{\bigg{(}}\int_t^T\!{\mathopen{\Big{(}}}\:\!\abs*{x(s) - \hat{x}(s)} + \abs*{\u(s) - \hat{\u}(s)}\:\!{\mathclose{\Big{)}}}\:\!ds\mathclose{\bigg{)}}}^{\!p\:\!}\:\!\! \\[1ex]
 &\lesssim \E{\abs*{x(T) - \hat{x}(T)}}^{\:\!p}\:\!\! + \E \int_t^T\!{\abs*{x(s) - \hat{x}(s)}}^{\:\!p}\:\!\!\:\!ds + \E \int_t^T\!{\abs*{\u(s) - \hat{\u}(s)}}^{\:\!p}\:\!ds \\[2.5ex]
 &\lesssim \E\:\!\!\sup_{s\in\mathopen{[}t,T\mathclose{]}}\:\!\!{\abs*{x(s) - \hat{x}(s)}}^{\:\!p}\:\!\! + \dist(\u(\:\!\bm{\cdot}\:\!),\hat{\u}(\:\!\bm{\cdot}\:\!)) \:\! .
\end{split}
\end{equation*}

In particular, the processes $\rounds{y(\:\!\bm{\cdot}\,;t),z(\:\!\bm{\cdot}\,;t)}$ and $\rounds{\hat{y}(\:\!\bm{\cdot}\,;t),\hat{z}(\:\!\bm{\cdot}\,;t)}$ no longer appear in the estimate. \textcolor{black}{Here, $\lesssim$\:\! indicates inequality up to a positive multiplicative constant, independent of the quantities involved, as we are not interested in specifying it further.} We also exploit the boundedness of $U$ to handle the term ${\abs*{\u(\:\!\bm{\cdot}\:\!) - \hat{\u}(\:\!\bm{\cdot}\:\!)}}^{\:\!p}$\:\!\!, once we observe that
\[
\E \int_t^T\!{\abs*{\u(s) - \hat{\u}(s)}}^{\:\!p}\:\!ds = \E \int_t^T\:\!\!\!\1_{\Set{\:\!\!\!\u(\:\!\bm{\cdot}\:\!) \:\! \neq \:\! \hat{\u}(\:\!\bm{\cdot}\:\!)\:\!\!\!}}(s,\omega) \:\! {\abs*{\u(s) - \hat{\u}(s)}}^{\:\!p}\:\!ds
\]
(see~Definition~\ref{def:dist}).

We now claim that
\begin{equation}
\label{eq:claim}
\E\:\!\!\sup_{s\in\mathopen{[}t,T\mathclose{]}}\:\!\!{\abs*{x(s) - \hat{x}(s)}}^{\:\!p}\:\!\! \;\! \lesssim \;\! \E \int_t^T\!{\abs*{\u(s) - \hat{\u}(s)}}^{\:\!p}\:\!ds
\end{equation} 
(that is, $\lesssim \dist(\u(\:\!\bm{\cdot}\:\!),\hat{\u}(\:\!\bm{\cdot}\:\!))$, as in the preceding step), after which the argument proceeds as outlined above, thereby concluding the proof.

Indeed, as in the previous estimates, for any $s\in\mathopen{[}t,T\mathclose{]}$,
\begin{multline*}
{\abs*{x(s) - \hat{x}(s)}}^{\:\!p}\:\!\! \:\! \lesssim \:\! \int_t^{\:\!s}\! {\abs*{\:\!b(r,x(r),\u(r)) - b(r,\hat{x}(r),\hat{\u}(r))}}^{\:\!p}\:\!\!\:\!dr \\
	+ {\abs*{\int_t^{\:\!s}\! {\mathopen{\big{(}}}\sigma(r,x(r),\u(r)) - \sigma(r,\hat{x}(r),\hat{\u}(r)){\mathclose{\big{)}}}\:\!dW(r)}}^{\:\!p} .
\end{multline*}
Hence, for any $\tau\in\mathopen{[}t,T\mathclose{]}$, applying the classical Burkholder\;\!--\:\!Davis\:\!--\:\!Gundy inequality over $\mathopen{[}t,T\mathclose{]}$ with exponent $r = p/2$, and then Jensen's inequality (as $p \geq 2$), we obtain
\begin{equation*}
\begin{split}
\hspace{0.001ex} & \E\:\!\!\sup_{s\in\mathopen{[}t,\tau\mathclose{]}}\:\!\!{\abs*{x(s) - \hat{x}(s)}}^{\:\!p}\:\!\! \\[0.75ex]
 &\lesssim \E \int_t^{\:\!\tau}\:\!\!{\abs*{b(r,x(r),\u(r)) - b(r,\hat{x}(r),\hat{\u}(r))}}^{\:\!p}\:\!\!\:\!dr \\[0.25ex]
 &+ \E\:\!\!\sup_{s\in\mathopen{[}t,\tau\mathclose{]}}{\abs*{\int_t^{\:\!s}\!{\mathopen{\big{(}}}\sigma(r,x(r),\u(r)) - \sigma(r,\hat{x}(r),\hat{\u}(r)){\mathclose{\big{)}}}\:\!dW(r)}}^{\:\!p}\:\!\! \\[0.75ex]
 &\lesssim \E \int_t^{\:\!\tau}\:\!\!{\abs*{b(r,x(r),\u(r)) - b(r,\hat{x}(r),\hat{\u}(r))}}^{\:\!p}\:\!\!\:\!dr \\[0.25ex]
 &+ \E\:\!\!{\mathopen{\bigg{(}}\int_t^{\:\!\tau}\!{{\mathopen{\big{(}}}\sigma(r,x(r),\u(r)) - \sigma(r,\hat{x}(r),\hat{\u}(r)){\mathclose{\big{)}}}}^{\:\!\!2}\:\!dr\mathclose{\bigg{)}}}^{\!p/2}\:\!\! \\[0.75ex]
 &\leq \E \int_t^{\:\!\tau}\:\!\!\mathopen{\Big{(}}\:\!{\abs*{b(r,x(r),\u(r)) - b(r,\hat{x}(r),\hat{\u}(r))}}^{\:\!p}\:\!\! + {\abs*{\sigma(r,x(r),\u(r)) - \sigma(r,\hat{x}(r),\hat{\u}(r))}}^{\:\!p}\:\!\mathclose{\Big{)}}\:\!dr \\[1ex]
 &\lesssim \E \int_t^{\:\!\tau}\!{\abs*{x(r) - \hat{x}(r)}}^{\:\!p}\:\!\!\:\!dr + \E \int_t^{\:\!\tau}\!{\abs*{\u(r) - \hat{\u}(r)}}^{\:\!p}\:\!dr \\[1ex]
 &\leq \int_t^{\:\!\tau}\!\E\:\!\!\sup_{s\in\mathopen{[}t,r\mathclose{]}}\:\!\!{\abs*{x(s) - \hat{x}(s)}}^{\:\!p}\:\!\!\:\!dr + \E \int_t^{\:\!\tau}\!{\abs*{\u(r) - \hat{\u}(r)}}^{\:\!p}\:\!dr \:\! .
\end{split}
\end{equation*}

Therefore, invoking Gr\"onwall's inequality in its integral form to the continuous function
\[
\tau \:\! \mapsto \:\! \E\:\!\!\sup_{s\in\mathopen{[}t,\tau\mathclose{]}}\:\!\!{\abs*{x(s) - \hat{x}(s)}}^{\:\!p}\:\!\!
\]
defined on $\mathopen{[}t,T\mathclose{]}$, we conclude that
\[
\E\:\!\!\sup_{s\in\mathopen{[}t,\tau\mathclose{]}}\:\!\!{\abs*{x(s) - \hat{x}(s)}}^{\:\!p}\:\!\! \;\! \lesssim \;\! \E \int_t^{\:\!\tau}\!{\abs*{\u(r) - \hat{\u}(r)}}^{\:\!p}\:\!dr \:\! ,
\]
which yields the desired estimate, namely~\eqref{eq:claim}, upon choosing $\tau = T$.




\end{appendices}





\end{document}